\documentclass[11pt]{article}
\usepackage{numbysec,amsfonts,amssymb,amsmath}
\hsize \textwidth
\advance \hsize by -\marginparwidth
\evensidemargin \oddsidemargin
\usepackage{amssymb}
\advance\hoffset by 5mm

\newtheorem{thm}[theorem]{Theorem}
\newtheorem{prop}[theorem]{Proposition}
\newtheorem{lemma}[theorem]{Lemma}
\newtheorem{remark}[theorem]{Remark}

\newtheorem{corollary}[theorem]{Corollary}
\newcommand{\Rm}{{\mathbb R}}
\newcommand{\pdr}[2]{\frac{\partial{#1}}{\partial{#2}}}

\newcommand{\eps}{\varepsilon}

\newcommand{\vu}{{\bf u}}
\newcommand{\vx}{{\bf x}}
\newcommand{\bv}{{\bf v}}
\newcommand{\bq}{{\bf q}}
\newcommand{\ve}{{\bf e}}

\newcommand{\commentout}[1]{}
\newcommand{\disp}{\displaystyle}
\numberbysection

\begin{document}
\title{Non-Planar Fronts in Boussinesq
Reactive Flows}
\author{
Henri Berestycki\footnote{EHESS, CAMS,
54 Boulevard Raspail, F - 75006 Paris, France; hb@ehess.fr}
\and Peter Constantin \footnote{Department of Mathematics, University of
Chicago, Chicago, IL 60637, USA; ryzhik@math.uchicago.edu}
\and Lenya Ryzhik\footnote{Department of Mathematics, University of
Chicago, Chicago, IL 60637, USA; const@math.uchicago.edu}}

\maketitle

\begin{abstract} We consider the reactive Boussinesq equations in
a slanted cylinder, with zero stress boundary conditions and arbitrary
Rayleigh number. We show that the equations have
non-planar traveling front solutions that propagate at a constant
speed. We also establish uniform upper bounds on the burning rate and the flow velocity
for general front-like initial data for the Cauchy problem.
\end{abstract}

\section{Introduction}

The existence of traveling fronts for reaction-diffusion equations and their stability has been extensively studied since the pioneering work of
Kolmogorov, Petrovskii and Piskunov \cite{KPP} and Fisher \cite{Fisher}.
A large number of results have been obtained during the last decade
on the generalization of the notion of a traveling front  to
reaction-diffusion-advection equations in a prescribed flow. These include non-planar
traveling fronts
in shear flows \cite{BLN,BLL,BN}, and pulsating traveling fronts
in periodic flows \cite{BH,X1,X2}, as well as results for monotonic
systems   in a unidirectional flow \cite{Volpert-1,Volpert-2,Volpert-3}.
One of the main qualitative effects of a flow is the speed-up of front
propagation due to front stretching. Various bounds have been obtained
for the speed of propagation of fronts in prescribed
flows \cite{ACVV1,ACVV2,ABP,BHN-1,CKOR,HPS,KR,KS,KRS,PX-91},
including variational principles for the front speed
\cite{BHN-1,BHN-2,F1,F2,Hamel,HPS}. The homogenization
limit in a periodic flow has also been studied \cite{MS}. Extensive
recent overviews can be found in \cite{B1,curved-fronts,X3}.

However, those results have been obtained under the
assumption that the flow is imposed from outside, and that it is not
affected by the evolution of the solution of the
reaction-diffusion-advection equation, that is, by the temperature
or concentration of the reactant. This is known as the constant
density approximation in the combustion literature. A first step
in the coupling of the temperature and fluid flow evolution is
via the Boussinesq approximation: the density mismatch is so small
that the density difference is accounted   by a buoyancy force
in the equation for an incompressible flow.
Recently a
number of works considered systems of a reaction-diffusion-advection
equation coupled to a flow equation of the
Boussinesq type. Global
existence and regularity of solutions in two dimensions was
studied in \cite{MX}. It has been shown that
non-planar convective traveling fronts may not exist in a vertical
cylinder if the Rayleigh number is too small while for large Rayleigh
numbers the planar fronts become unstable \cite{CKR,TV1,TV2}.
Moreover, there exists a bifurcation at a critical value $\rho_c>0$ --
non-trivial convective fronts may exist for the Rayleigh numbers
close to $\rho_c$ \cite{TV1,TV2}. Numerical computations \cite{VR} show that
non-planar convective fronts exist and are stable for a large range of
Rayleigh numbers $\rho>\rho_c$. The fingering instability in this
regime was investigated in \cite{deWit}.

One of the difficulties in the analysis of the Boussinesq problem
at large Rayleigh numbers in a vertical cylinder is the presence
of unstable planar fronts that make uniform lower bounds on the
front speed quite difficult. However, it has been observed in
\cite{Volpert-private} that such planar fronts cannot exist in a
horizontal cylinder. One of the main results of
\cite{Volpert-private} is that non-planar fronts in a horizontal
cylinder exist for small Rayleigh numbers. A purpose of the present
paper is to extend this result to all positive Rayleigh numbers;
we use an approach that is different from \cite{Volpert-private} and
is based on the a priori bounds developed in \cite{CKR}.

The reactive Boussinesq equations for the temperature $T$ and flow
$\vu$ have the dimensional form
\begin{eqnarray}\label{bouss-dim}
T_t+\vu\cdot\nabla T=\kappa\Delta T+\frac{v_0^2}{\kappa}f(T)\\
\vu_t+\vu\cdot\nabla\vu-\nu\Delta\vu+\nabla p= gT{\bf e}_z\nonumber\\
\nabla\cdot \vu=0.\nonumber
\end{eqnarray}
Here ${\bf e}_z$ is the unit vector in the vertical direction, $g$ is
the strength of gravity, the speed $v_0$ is proportional to the traveling front speed in
the absence of
gravity, $\kappa$ is the thermal diffusivity and $\nu$ is the fluid viscosity.
The temperature is normalized so that $0\le T\le 1$.
The nonlinearity $f(T)$ is assumed to be a Lipschitz function of the ignition type
\begin{equation}\label{f-comb}
\hbox{$f(T)=0$ for $0\le T\le\theta_0$ with $\theta_0>0$,
$f(T)>0$ for $T\in(\theta_0,1)$ and $f(1)=0$.}
\end{equation}

We consider the equations (\ref{bouss-dim}) in a slanted two-dimensional cylinder $x\in
\Rm$, $\alpha x\le z\le \alpha x+H$ with a finite slope
$\alpha<\infty$.  It is convenient to rotate the cylinder in order to
make it horizontal to simplify the notation. Then (\ref{bouss-dim})
becomes
\begin{eqnarray}\label{bouss-dim-2}
T_t+\vu\cdot\nabla T=\kappa\Delta T+\frac{v_0^2}{\kappa}f(T)\\
\vu_t+\vu\cdot\nabla\vu-\nu\Delta\vu+\nabla p=
gT\hat\ve \nonumber\\
\nabla\cdot \vu=0,\nonumber
\end{eqnarray}
where $\vu$ is the flow velocity measured relative to the new
coordinate system.  The gravity on the right points in a direction
$\hat\ve$ that is   non-parallel to the $x$-axis, as the original cylinder was assumed to
be non-vertical ($\alpha<\infty$). The new rotated problem is posed in
a cylinder $D=\Rm_x\times[0,L]_z$, $L=H/\sqrt{1+\alpha^2}$.  The boundary
conditions for the temperature $T$ are set to be front-like:
\begin{equation}\label{bc-T}
\hbox{$T\to 1$ as
$x\to-\infty$, $T\to 0$ as $x\to +\infty$, $\disp\pdr{T}{z}=0$
at $z=0,L$.}
\end{equation}
The flow $\vu=(v,w)$ satisfies the no stress boundary conditions:
\begin{equation}\label{bc-u}
\hbox{$\vu,\omega\to 0$ as $x\to\pm\infty$ and $w,\omega=0$ at $z=0,L$.}
\end{equation}
Here $\omega=w_x-v_z$ is the flow vorticity so that
\[
\Delta v=-\omega_z,~~\Delta w=\omega_x.
\]
In order to pass to the non-dimensional variables we introduce the laminar
front width $\delta=\kappa/v_0$ and reaction time $t_c=\kappa/v_0^2$
and rescale the space and time variables: $\vx_{new}=\vx_{old}/\delta$
and $t_{new}=t_{old}/t_c$. We also rescale the flow $\vu_{new}=\vu_{old}/v_0$.
Then the Boussinesq equations become
\begin{eqnarray}\label{bouss-nondim}
T_t+\vu\cdot\nabla T=\Delta T+f(T)\\
\vu_t+\vu\cdot\nabla\vu-\sigma\Delta\vu+\nabla p=
\rho T\hat\ve\nonumber\\
\nabla\cdot \vu=0,\nonumber
\end{eqnarray}
where $\sigma=\nu/\kappa$ is the Prandtl number and
$\rho= {g\delta^3}/{\kappa^2}$ is
the Rayleigh number. The problem is now posed in the strip
$D=\Rm_x\times[0,\lambda]_z$, $\lambda=L/\delta$, with the boundary
conditions that come from (\ref{bc-T}) and (\ref{bc-u}).

The traveling front solutions of (\ref{bouss-nondim}) are solutions of the form
$T(x-ct,z)$, $\vu(x-ct,z)$ with the speed $c$ to be determined. They satisfy
\begin{eqnarray}\label{bouss-travel}
-cT_x+\vu\cdot\nabla T=\Delta T+f(T)\\
-c\vu_x+\vu\cdot\nabla\vu-\sigma\Delta\vu+\nabla p=
\rho T\hat\ve\nonumber\\
\nabla\cdot \vu=0,\nonumber
\end{eqnarray}
with the boundary conditions
\begin{equation}\label{bc-T-trav}
\hbox{$T\to T_-$ as $x\to-\infty$, $T\to 0$ as $x\to +\infty$,
$\disp\pdr{T}{z}=0$ at $z=0,\lambda$}
\end{equation}
and
\begin{equation}\label{bc-u-trav}
\hbox{$w,\omega=0$ at $z=0,\lambda$.}
\end{equation}
Here $T_-$ is a constant that is not a priori prescribed.
We recall that, as has been  observed in \cite{Volpert-private},
if the direction of gravity $\hat\ve$ is not parallel to the $x$-axis,
any traveling front solution of (\ref{bouss-travel})
must be non-planar, that is, it must depend on both variables $x$ and
$z$. This is the main difference between the cases of a vertical and
slanted cylinder: planar fronts exist in the former case but not in
the latter.  Our main
result is the following theorem.
\begin{thm}\label{thm-main}
Let the nonlinearity $f(T)$ be of the ignition type
(\ref{f-comb}). Then a traveling front solution
$(c,T,\vu)$ of (\ref{bouss-travel}) exists such that it is
non-planar: $T_z\not\equiv 0$, the flow
$\vu\not\equiv 0$ and the reaction rate $f(T)\not\equiv 0$.
Moreover, the solution satisfies the following properties:
$c>0$, $T\in C^{2,\alpha}(D)$, $\nabla T\in L^2(D)$, $\vu\in H^1(D)\cap C^{2,\alpha}(D)$.
If we assume in addition that
\begin{equation}\label{f-technical}
\hbox{$f(T)\le(T-\theta_0)_+^2/\lambda^2$,}
\end{equation}
then the left limit $T_-=1$.
\end{thm}
The assumption (\ref{f-technical}) is of technical nature. It
does not involve the Rayleigh number $\rho$,
it is rather a restriction on the channel width $\lambda$. We do
not address the question of the uniqueness of the traveling front
speed or profile in this paper -- this problem requires an additional
study. Our results can be generalized to the no-slip boundary
conditions $\vu=0$ on $\partial D$ at the expense of a more
technical proof -- we leave this problem for a future
publication.

The general idea of the proof is as follows. We first consider the problem
(\ref{bouss-travel}) on a finite domain
$D_a=[-a,a]_x\times[0,\lambda]$.  Solutions $(T_a^c,\vu_a^c)$ of the
restricted problem exist for all $c\in\Rm$. We normalize them by the
requirement that
\begin{equation}\label{normal}
\disp\max_{x\ge 0}T_a^c(x,z)=\theta_0.
\end{equation}
This imposes a restriction on the speed $c$.  In order to show that
there exists a speed $c_a$ so that (\ref{normal}) holds we first obtain some a
priori bounds on $c$, $T$ and $\vu$ under the condition
(\ref{normal}). Then we use the Leray-Schauder topological degree
theory and the above a priori bounds to show that $c_a$ exists. The a
priori bounds allow us to pass to the limit $a\to\infty$. Finally we
show that the right limit of $T$ as $x\to+\infty$ is equal to zero,
and that the left limit is equal to one under the
additional assumption on $f(T)$ in Theorem \ref{thm-main}. This general
strategy is similar to that in the proof of existence of traveling fronts in
a prescribed decoupled flow, as in, for example, \cite{BLN,BN}.
The main difficulty and novelty are in the a priori bounds for the
solution of the coupled problem in a bounded domain.

Our second result shows that the solution of the Cauchy problem
for (\ref{bouss-nondim}) propagate with a finite speed and
that this speed is close to the speed of the laminar front $c_0$ when
the Rayleigh number is small. Recall that there exists a unique
speed $c_0$ so that a traveling front solution of
\[
-c_0\Phi_x=\Phi_{xx}+\Phi(U),~~\Phi(-\infty)=1,~\Phi(+\infty)=0
\]
exists.

In order to make this precise we define the bulk burning rate $\bar V(t)$, the
Nusselt
number $\bar N(t)$ and the average horizontal flow $\bar U(t)$ by
\begin{eqnarray}\label{intro-4-def}
&&\bar V(t)=\frac{1}{t}\int_0^t\int V(s)ds,~~V(t)=\int f(T)\frac{dxdz}{\lambda},\\
&&\bar N(t)=\frac{1}{t}\int_0^tN(s)ds,~~N(t)=\int |\nabla
T|^2\frac{dxdz}{\lambda},\\
&&\bar U(t)=\frac{1}{t}\int_0^t\|v(s)\|_\infty ds.
\end{eqnarray}
The following theorem provides uniform bounds on these bulk
quantities. It also shows that the coupled problem (\ref{bouss-nondim})
is in a sense a "regular perturbation" of the single
reaction-diffusion equation with $\rho=0$.
\begin{thm}\label{intro-thm-bds} Assume that
there exists $R$ so that $T_0(x,z)=0$ for $x>R$ and $T_0(x,z)=1$ for
$x<-R$ and that the initial vorticity $\omega_0\in L^2(D)$. There exists a constant $C>0$
so that
under the above assumptions on the initial data $T_0$, $\vu_0$ we
have the following bounds
\begin{eqnarray}\label{intro-4-bds}
&&c_0-C[\rho+\rho^2]+o(1)\le\bar V(t)\le c_0+C[\rho+\rho^2]+o(1)\\
&&\bar N(t)\le \left[{C\rho}+\sqrt{\frac{c_0}{2}+C^2\rho^2}\right]^2+o(1)\nonumber\\
&&\bar U(t)\le C\rho[1+\rho]+o(1)\nonumber
\end{eqnarray}
as $t\to +\infty$.
\end{thm}
This theorem may be interpreted as a
stability result for a perturbation of a homogeneous reaction-diffusion equation
by the buoyancy coupling. The proof is based on the construction
of super- and sub- solutions, and a bound on the decay of the
solutions of advection-diffusion equations that is uniform in the
advection flow.

The third result of this paper deals with the Boussinesq system in a narrow
domain. It has been shown in \cite{CKR} that if a vertical strip
is sufficiently narrow and gravity is sufficiently weak then
solutions of the Cauchy data become planar as $t\to +\infty$. The
following theorem generalizes this result to inclined cylinders.
\begin{thm}\label{intro-thm-narrow} Let $\hat\ve=(e_1,e_2)$ be the unit
vector in the direction of gravity and let $\rho_j=\rho e_j$,
$j=1,2$ and let the initial data $(T_0,\vu_0)$ be as in
Theorem \ref{intro-thm-bds}. There exist two constants $\lambda_0$ and $\rho_0$ so
that if the domain is sufficiently narrow: $\lambda\le\lambda_0$
and gravity is sufficiently small: $\rho\le\rho_0$ then the
burning rate is bounded by
\begin{equation}\label{intr-5-barv}
\bar V(t)\le c_0+C\rho_2+o(1)\hbox{ as $t\to +\infty$.}
\end{equation}
Moreover, the front is nearly planar in the sense that
\begin{equation}\label{intro-5-Tx}
\bar N_z(t)=\frac{1}{t}\int_0^t\|T_z(s)\|_2^2ds\le C\rho_2^2+o(1)
\hbox{ as $t\to +\infty$.}
\end{equation}
\end{thm}
The main observation of this theorem is that only the gravity strength in
the direction perpendicular to the strip enters in the upper
bounds (\ref{5-barv}) and (\ref{5-Tx}).

The paper is organized as follows: Theorem \ref{thm-main} is proved
in Sections \ref{sec:finite} and \ref{sec:limit}. Theorems
\ref{intro-thm-bds} and \ref{intro-thm-narrow} are proved in Section \ref{sec:bds}.

{\bf Acknowledgment.} We thank Vitaly Volpert for explaining to us
the results of \cite{Volpert-private} prior to its publication. We
also thank Marta Lewicka for a careful reading of the preliminary
version of the manuscript.
This research was supported in part by the
ASCI Flash center at the University of Chicago under DOE contract B341495.
PC was partially supported by the NSF grant DMS-0202531, LR by NSF grant
DMS-0203537, ONR grant N00014-02-1-0089 and an Alfred P. Sloan Fellowship.

\section{The finite domain problem}\label{sec:finite}

We consider in this section the approximating problem
\begin{eqnarray}\label{bouss-travel-finite}
-cT_x+\vu\cdot\nabla T=\Delta T+f(T)\\
-c\vu_x+\vu\cdot\nabla\vu-\sigma\Delta\vu+\nabla p=\rho T\hat{\bf e}\nonumber\\
\nabla\cdot \vu=0,\nonumber
\end{eqnarray}
in a finite domain $D_a=[-a,a]_x\times[0,\lambda]_y$, $a>0$, with the
boundary conditions
\begin{equation}\label{bc-T-trav-fin}
\hbox{$T(-a,z)=1$,  $T(a,z)=0$,
$\disp\pdr{T}{z}=0$ at $z=0,\lambda$}
\end{equation}
and
\begin{equation}\label{bc-u-trav-fin}
\hbox{$w=0$, $\omega=0$ at $z=0,\lambda$ and $v(\pm a,z)=\omega(\pm a,z)=0$
at $x=\pm a$.}
\end{equation}
One can show with the techniques of the present section that a solution
$T_a$, $\vu_a$ of (\ref{bouss-travel-finite}) in $D_a$
with the boundary conditions (\ref{bc-T-trav-fin}) and
(\ref{bc-u-trav-fin}) exists for all $c\in\Rm$. However, given an
arbitrary $c$ there is no way to control the limit of $T_a$ and
$\vu_a$ as $a\to \infty$. Hence, following the standard procedure,
we impose an additional constraint (\ref{normal}).
This ensures that the non-trivial part of the
solution does not escape to infinity when we pass to the limit
$a\to\infty$.
\begin{prop}\label{lem-1/2}
There exists a speed $c_a\in\Rm$ so that there exists a solution
$(T_a,\vu_a)$ of (\ref{bouss-travel-finite}) in $D_a$
with the boundary conditions (\ref{bc-T-trav-fin}) and
(\ref{bc-u-trav-fin}) such that
\begin{equation}\label{max}
\max_{x\ge0,z\in[0,\lambda]}T_a(x,z)=\theta_0.
\end{equation}
We denote the corresponding solution as $(c_a,T_a,\vu_a)$.  Moreover,
there exists $a_0>0$ and a constant $C>0$ that
is independent of $a$, so that we have for all $a>a_0$
\begin{equation}\label{c-bd-a}
|c_a|\le C,
\end{equation}
and
\begin{equation}\label{uT-bd-a}
\int_{D_a}|\nabla T_a|^2dxdz+
\int_{D_a}|\nabla \vu_a|^2dxdz+
\|\vu_a\|_\infty\le C.
\end{equation}
Moreover, the uniform H\"older estimates hold:
there exists $a_0>0$ and a constant $C>0$ independent
of $a$ so that we have for
all $a>a_0$
\begin{equation}\label{gradT-infty-holder}
\|\omega_a\|_{C^{1,\alpha}(D_a)}+\|\vu_a\|_{C^{1,\alpha}(D_a)}+
\|T_a\|_{C^{1,\alpha}(D_a)}\le C
\end{equation}
provided that $0<\alpha<1$.
\end{prop}
{\bf Proof.}  The proof consists of two parts. First, we introduce a
family of problems depending on a parameter $\tau\in[0,1]$ so that at
$\tau=0$ we have a simple linear problem without advection or coupling
and at $\tau=1$ we have the full problem (\ref{bouss-travel-finite})
with the correct boundary conditions.  The normalization condition
(\ref{max}) is imposed for all $\tau\in[0,1]$.  We obtain the a priori
bounds as in (\ref{c-bd-a}), (\ref{uT-bd-a}) and
(\ref{gradT-infty-holder}) for such solutions that are uniform in
$\tau\in[0,1]$. In the second step we use the a priori bounds, the
Leray-Schauder topological degree argument and the information on the
linear problem at $\tau=0$ to show that solutions of the nonlinear
coupled problem at $\tau=1$ exist.  We drop the subscript $a$
throughout the proof to make the notation less cumbersome.

{\bf Step 1. A priori bounds for solutions.} Let us first define a
one-parameter (homotopy) family of finite domain Boussinesq problems
in the vorticity formulation
\begin{eqnarray}\label{1.1.tau-2}
&&-c^\tau T_x^\tau+\tau\vu^\tau\cdot\nabla T^\tau=\Delta T^\tau+\tau f(T^\tau)\\
&&-c^\tau\omega_x^\tau+\vu^\tau\cdot\nabla\omega^\tau-\sigma\Delta\omega^\tau=\tau\rho\hat{\bf
e}\cdot\nabla^\perp T:= \rho\tau[e_2 T_x^\tau
-e_1T_z^\tau]\nonumber\\
&&\omega^\tau=w_x^\tau-v_z^\tau,~~\nabla\cdot\vu^\tau=0.\nonumber
\end{eqnarray}
As mentioned above, $\tau$ is the homotopy parameter: $\tau\in[0,1]$,
with $\tau=0$ corresponding to the linear problem, and $\tau=1$ to the
full problem (\ref{bouss-travel-finite})-(\ref{bc-u-trav-fin}).
The problem (\ref{1.1.tau-2}) is posed in
$D_a$ with the same boundary conditions
\begin{equation}\label{bc-lat.tau-2}
\pdr{T^\tau}{z}=0,~~w^\tau=\omega^\tau=0~~\hbox{for $z=0,\lambda$}
\end{equation}
and
\begin{equation}\label{bc-lat-a.tau-2}
T^\tau(-a,z)=1,~~T^\tau(a,z)=0,~~
v^\tau(\pm a,z)=\omega^\tau=(\pm a,z)=0~~\hbox{for $x=\pm a$},
\end{equation}
as (\ref{bouss-travel-finite}).
We also require that
\begin{equation}\label{T1/2.tau}
\max_{x\ge 0,z} T^\tau(x,z)=\theta_0
\end{equation}
and obtain a priori bounds on $c^\tau$, $T^\tau$ and
$\omega^\tau$. We drop the superscript $\tau$ below wherever it
causes no confusion. The general plan is as follows. First, we
bound the speed $c$ above and below by a linear function of
$\|v\|_\infty$ in Lemma \ref{lem-cu}.  Next we bound
$\|\vu\|_\infty$ from above by a linear function of $\|\nabla
T\|_2$ in Lemma \ref{lem-uTbd}. The other direction, a bound on
$\|\nabla T\|_2^2$ in terms of a linear function of
$\|\vu\|_\infty$ is established in Lemmas \ref{lem-Tu-bd} and
\ref{lem-Tx}. Since the latter bound is quadratic in $\|\nabla
T\|_2$, the last estimates allow to obtain a uniform bound on this
quantity, from which all other a priori bounds follow in a fairly
straightforward manner: see Corollary \ref{cor-bds} and Lemma
\ref{lem-holder}.

We begin with a lemma that bounds the speed $c$ in terms of the
horizontal flow velocity  $\|v\|_{L^\infty(D_a)}$.
\begin{lemma}\label{lem-cu}
Let $(c,T,\vu)$ satisfy (\ref{1.1.tau-2})-(\ref{bc-lat-a.tau-2}) with
the normalization (\ref{T1/2.tau}) and let $\vu=(v,w)$. Then there
exists $a_0>0$ so that for all $a\ge a_0$ we have
\begin{equation}\label{ctau-bds}
-1-\tau\|v\|_\infty\le
c\le 1+M\tau+\tau\|v\|_\infty.
\end{equation}
\end{lemma}
{\bf Proof.} First, we observe that the function
$\psi_A(x)=Ae^{-\alpha(x+a)}$ is a super-solution for the
reaction-diffusion-advection equation with the flow $\vu$ fixed if
$A>1$ and
\begin{equation}\label{calpha}
c\ge \alpha+\frac{M\tau}{\alpha}+\tau\|v\|_\infty,
\end{equation}
that is,
\begin{equation}\label{super}
-c\pdr{\psi_A}{x}+\tau\vu\cdot\nabla\psi_A\ge\Delta\psi_A+\tau f(\psi_A),
\end{equation}
provided that (\ref{calpha}) holds with
\[
M=\sup_{0\le T\le 1}\frac{f(T)}{T}.
\]
Furthermore, we have
\begin{equation}\label{super-endpts}
T(-a,z)=1<A=\psi_A(-a),~~T(a,z)=0<\psi_A(a)
\end{equation}
at the two ends of the domain $D_a$. We now show that this together with (\ref{super})
implies that
\begin{equation}\label{T-psiA}
T(x,z)\le \psi_A(x)
\end{equation}
for all $(x,z)\in D_a$ and $A>1$. Indeed, consider the family of
functions $\psi_A(x)$. Then all $\psi_A$ are super-solutions in
the sense that the inequality (\ref{super}) holds. Moreover, as
the maximum principle implies that $0\le T\le 1$, for
$A>5e^{2\alpha a}$ sufficiently large we have $\psi_A(x)>5>T(x,z)$
for all $(x,z)\in D_a$. We define
\[
A_0=\inf\left\{A\in\Rm:~~\psi_A(x)\ge T(x,z)\hbox{ for all $(x,z)\in D_a
$}\right\}.
\]
The previous argument implies that $A_0$ is finite, $A_0\le
5e^{2\alpha a}$ and, moreover, clearly $A_0>0$. Observe that since
the domain $D_a$ is compact, we should have $\psi_{A_0}(x)\ge
T(x,z)$ -- otherwise this inequality would be violated for $A$
slightly larger than $A_0$ at some point in $D_a$. Moreover, the
equation $\psi_{A_0}(x)=T(x,z)$ should have a solution. We claim
that $A_0=1$. Indeed, otherwise the point $(x_0,z_0)$ that solves
$\psi_{A_0}(x_0)=T(x_0,z_0)$ cannot be at the boundary of $D_a$
because of the boundary conditions on the function $T$. Hence this
point has to lie in  the interior of $D_a$.  The continuity of
$\psi_A(x)$ with respect to $A$ implies that the graphs of
$\psi_{A_0}(x)$ and $T(x,z)$ are tangent at $(x_0,z_0)$. Then the
strong maximum principle implies that $\psi_{A_0}(x)\equiv T(x,z)$
which is a contradiction, as they differ on the boundary. Hence we
conclude that $A_0=1$ and thus (\ref{T-psiA}) holds for all $A>1$
and thus for $A=1$, so that
\begin{equation}\label{T-psi}
T(x,z)\le e^{-\alpha(x+a)}.
\end{equation}

However, the existence of such a super-solution contradicts
the normalization condition (\ref{T1/2.tau}) if $\alpha\ge \ln (\theta_0^{-1})/a$ because
(\ref{T1/2.tau}) implies that there exists $z_0$ so that
$T(0,z_0)=\theta_0$. Therefore,  the  existence of a solution $T$ that
satisfies (\ref{T1/2.tau}) implies
\begin{eqnarray}
&&c\le \inf_{\alpha\ge \ln (\theta_0^{-1})/a}
\left(\alpha+\frac{M\tau}{\alpha}\right)+\tau\|v\|_\infty\le
1+M\tau +\tau\|v\|_\infty
\label{c-tau-up}
\end{eqnarray}
provided that $a\ge\ln(1/\theta_0)$.
This proves the upper bound in (\ref{ctau-bds}). In order to prove the
lower bound we observe that the function $\phi=1-e^{\alpha (x-a)}$ is
a sub-solution for $T$ with the flow $\vu$ fixed if
\begin{equation}\label{ctau-low-0}
c\le -\alpha-\tau\|v\|_\infty.
\end{equation}
That is, if (\ref{ctau-low-0}) holds, then $T(x,z)\ge 1-e^{\alpha
(x-a)}$.  This is shown in a way similar to the proof of (\ref{T-psi})
under the assumption (\ref{calpha}) above.  However, $\phi(0)=1-e^{-\alpha
a}>\theta_0$ for
\begin{equation}\label{aln2}
a>\frac{\ln ((1-\theta_0)^{-1})}{\alpha}.
\end{equation}
This implies that $\max_{x\ge 0}T(x,z)\ge \phi(0)>\theta_0$
provided that both (\ref{ctau-low-0}) and (\ref{aln2}) hold. Hence in
order for (\ref{T1/2.tau}) to be possible we need
\begin{equation}\label{ctau-low}
c\ge \sup_{\alpha>\frac{\ln
((1-\theta_0)^{-1})}{a}}\left[-\alpha-\tau\|v\|_\infty\right]\ge
-1-\tau\|v\|_\infty
\end{equation}
provided that $a\ge \ln((1-\theta_0)^{-1})$.
This is the lower bound in (\ref{ctau-bds}) and the proof of Lemma
\ref{lem-cu} is complete. $\Box$

Next we establish a bound on $\|\vu\|_{L^\infty(D_a)}$ and
$\|\omega\|_{L^\infty(D_a)}$ in terms of $\|\nabla T\|_{L^2(D_a)}$.
These bounds are all obtained from the following type of estimates.
\begin{lemma}\label{lem-laplace}
Let $S_a=[-a,a]_x\times\Omega_y$ be a finite cylinder with a smooth
bounded cross-section $\Omega\in\Rm^d$, $d=1,2$. Let $\phi$ be a
function that satisfies either of the following three conditions: (i)
$\phi(x,y)=0$ on the whole boundary $\partial S_a$, (ii) $\phi(x,y)=0$
for $y\in\partial\Omega$ and $\pdr{\phi(x,y)}{x}=0$ for $x=-a,a$, or
(iii) $\pdr{\phi(x,y)}{n}=0$ for $y\in\partial\Omega$, and $\phi(x,y)=0$
for $x=-a,a$.  Then there exists a constant $C$ that depends only on
the domain $\Omega$, but not on the cylinder length $a$, so that we
have
\begin{equation}\label{linfty-laplace}
\|\phi\|_{L^\infty(S_a)}\le
C\left[\|\Delta\phi\|_{L^2(S_a)}+\|\phi\|_{L^2(S_a)}\right].
\end{equation}
\end{lemma}
{\bf Proof.} Let $Q$ be any cylinder of the form
$[x_0,x_0+1]\times\Omega\subset S_a$ with $-a\le x_0\le a-1$.  The
standard interior elliptic estimates up to the boundary \cite{GT} can
be applied to $Q$ in all the three cases (i)-(iii). The corners at
$x=\pm a$ are not an obstacle. Indeed, both in the case of the
Dirichlet and Neumann boundary conditions prescribed on the lines
$x=\pm a$, one can extend the solution to a larger cylinder
$[-a-1,a+1]\times\Omega$ by reflecting the solution across the line
$x=\pm a$, either in the even or odd way, respectively. Hence the
usual elliptic estimates up to the boundary can be applied to all such
cylinders $Q$ to obtain
\begin{equation}\label{h2-est}
\|\phi\|_{H^2(Q)}\le
C\left[\|\Delta\phi\|_{L^2(S_a)}+\|\phi\|_{L^2(S_a)}\right]
\end{equation}
in all three cases (i)-(iii). Then the Sobolev embedding theorem in
dimensions $d=2,3$ implies that
\[
\|\phi\|_{L^\infty(Q)}\le C\|\phi\|_{H^2(Q)}\le
C\left[\|\Delta\phi\|_{L^2(Q)}+\|\phi\|_{L^2(Q)}\right] \le
C\left[\|\Delta\phi\|_{L^2(S_a)}+\|\phi\|_{L^2(S_a)}\right]
\]
with the constant $C$ that depends only on the domain $\Omega$. $\Box$

This lemma can be easily extended to higher dimensions using the
appropriate Sobolev embeddings. It implies immediately the following
bounds on $\|\vu\|_\infty$ and $\|\omega\|_\infty$ in terms of
$\|\nabla T\|_{L^2(D_a)}$.
\begin{lemma}\label{lem-uTbd}
Let $(c,T,\vu)$ satisfy
(\ref{1.1.tau-2})-(\ref{bc-lat-a.tau-2}) with the normalization
(\ref{T1/2.tau}). There exists $a_0>0$ and a constant $C>0$ so that,
for all $a>a_0$, it holds that
\begin{equation}\label{uT-bd}
\|\vu\|_{L^\infty(D_a)}\le C\|\nabla T\|_{L^2(D_a)}.
\end{equation}
and
\begin{equation}\label{omega-infty-bd}
\|\omega\|_{L^\infty(D_a)}\le C\|\nabla T\|_{L^2(D_a)}
\left[1+\|\nabla T\|_{L^2(D_a)}\right].
\end{equation}
Moreover, $\nabla\vu$ satisfies the same bound:
\begin{equation}\label{nabla-vu-infty-bd}
\|\nabla\vu\|_{L^\infty(D_a)}\le C\|\nabla
T\|_{L^2(D_a)}\left[1+\|\nabla T\|_{L^2(D_a)}\right].
\end{equation}
\end{lemma}
{\bf Proof.}
We use the vorticity equation
\begin{equation}\label{vort-eq}
-c\omega_x+\vu\cdot\nabla\omega-
\sigma\Delta\omega=\rho\tau(\hat\ve\cdot\nabla^\perp T),
~~\hbox{$\omega=0$ on $\partial D_a$.}
\end{equation}
Case (i) of Lemma \ref{lem-laplace} implies that
\begin{equation}\label{vortinfty}
\|\omega\|_{L^\infty(D_a)}\le
C\left[\|\nabla T\|_{L^2(D_a)}+
(|c|+\|\vu\|_\infty)\|\nabla\omega\|_{L^2(D_a)}
+\|\omega\|_{L^2(D_a)}\right].
\end{equation}
Here the constant $C$ depends only on $\rho$ and $\lambda$.  Note
that multiplying the vorticity equation by $\omega$ and integrating by
parts, using the boundary conditions we obtain
\[
 \int_{D_a}|\nabla\omega|^2dxdz=
\tau\rho\int\left(\hat\ve\cdot\nabla^\perp T\right)\omega dxdz\le  \tau
\rho\|\nabla T\|_2\|\omega\|_2.
\]
The Dirichlet boundary conditions for $\omega$ imply that the Poincar\'e
inequality applies to $\omega$ so that $\|\omega\|_{L^2(D_a)}\le
(\lambda/\pi)\|\nabla\omega\|_{L^2(D_a)}$.  Hence we obtain
\begin{equation}\label{nablaomega-T}
\|\nabla\omega\|_2\le \frac{\lambda}{\pi}\tau\rho\|\nabla T\|_2
\end{equation}
and thus
\begin{equation}\label{omega-T}
\|\omega\|_{L^2(D_a)}\le C\|\nabla T\|_{L^2(D_a)},
\end{equation}
with the constant $C$ independent of the cylinder length $a$. This,
together with (\ref{vortinfty}) and the bound (\ref{ctau-bds})
on the speed $c$ implies (\ref{omega-infty-bd}), provided that we show
(\ref{uT-bd}).

We now prove (\ref{uT-bd}).
The horizontal flow component $v$ satisfies the
Poisson equation
\begin{equation}\label{u-poisson}
\Delta v=-\omega_z,~~v(\pm a,z)=0,~~\pdr{v}{z}=0,~~\hbox{at
$z=0,\lambda$.}
\end{equation}
The boundary conditions at $z=0,\lambda$ are obtained from
$v_z=w_x-\omega=0$ as follows from (\ref{bc-lat.tau-2}). The third case
(iii) of Lemma \ref{lem-laplace} implies that
\begin{equation}\label{uinfty}
\|v\|_{L^\infty(D_a)}\le C[\|\nabla \omega\|_{L^2(D_a)}+\|v\|_{L^2(D_a)}].
\end{equation}
The first term in the right side is bounded by (\ref{nablaomega-T}). In
order to bound the second we multiply (\ref{u-poisson}) by $v$ and integrate
to obtain, using the boundary conditions and (\ref{nablaomega-T})
\begin{equation}\label{v-nabla}
\int_{D_a}|\nabla v|^2dxdz=
\int_{D_a} \omega_z(x,z)v(x,z)dxdz\le\|\omega_z\|_2\|v\|_2\le
C\|\nabla T\|_2\|v\|_2.
\end{equation}
Now, observe that (\ref{u-poisson}), the Neumann boundary conditions for
$v$ and the Dirichlet boundary condition for $\omega$ at $z=0,\lambda$ imply
that
\[
\frac{d^2}{dx^2}\int v(x,z)dz=0.
\]
It follows then from the Dirichlet boundary conditions for $v$ at $x=\pm a$
that
\begin{equation}\label{v-meanzero}
\int_0^\lambda v(x,z)dz=0
\end{equation}
for all $x$. One may alternatively deduce (\ref{v-meanzero}) from
incompressibility of the flow $\vu$ and the boundary conditions.
Therefore, it follows from the Poincar\'e inequality that
$\|v\|_{L^2(D_a)}\le(\lambda/2\pi)\|\nabla v\|_{L^2(D_a)}$. Thus,
(\ref{v-nabla}) implies that both $\|\nabla v\|_{L^2(D_a)}\le C\|\nabla
T\|_{L^2(D_a)}$ and $\|v\|_{L^2(D_a)}\le C\|\nabla
T\|_{L^2(D_a)}$ with a constant independent of $a$. Hence,
(\ref{uinfty}) implies (\ref{uT-bd}) for the horizontal flow component.


The vertical flow component satisfies
\begin{equation}\label{vert-eq}
\Delta w=\omega_x,~~w(x,0)=w(x,\lambda)=0,~~\pdr{w}{x}(\pm a,z)=0.
\end{equation}
The Neumann boundary condition at $x=\pm a$ is deduced from the relation
$w_x=\omega+v_z$
and the Dirichlet boundary conditions for $v$ and $\omega$ at $x=\pm a$.
The case (ii) in Lemma \ref{lem-laplace} implies that
\begin{equation}\label{winfty}
\|w\|_{L^\infty(D_a)}\le C[\|\nabla \omega\|_{L^2(D_a)}+\|w\|_{L^2(D_a)}].
\end{equation}
As before, we use (\ref{nablaomega-T}) to bound the first term in the
right side. In order to bound the second we multiply (\ref{vert-eq})
by $w$ and integrate, using the boundary conditions and
(\ref{nablaomega-T}) again, to obtain that
\begin{equation}\label{w-nabla}
\int_{D_a}|\nabla w|^2dxdz=-
\int_{D_a} \omega_x(x,z)w(x,z)dxdz\le\|\omega_x\|_2\|w\|_2\le
C\|\nabla T\|_2\|w\|_2.
\end{equation}
The Dirichlet boundary conditions for $w$ at $z=0,\lambda$ imply that
$\|w\|_{L^2(D_a)}\le\lambda/\pi\|\nabla w\|_{L^2(D_a)}$. Thus,
(\ref{w-nabla}) implies that
\begin{equation}\label{nablaw-T}
\|\nabla w\|_{L^2(D_a)}\le C\|\nabla
T\|_{L^2(D_a)},
\end{equation}
and hence $\|w\|_{L^2(D_a)}\le C\|\nabla T\|_{L^2(D_a)}$ with a
constant independent of $a$. Therefore, now (\ref{winfty}) implies
(\ref{uT-bd}) for the vertical flow component. Thus, the proof of
(\ref{uT-bd}) is complete. We recall that then (\ref{omega-infty-bd})
follows as well, as explained in the paragraph below
(\ref{omega-T}).

In order to complete the proof of Lemma \ref{lem-uTbd} it remains to
bound the derivatives of $\vu$. First, we observe that the function
$\psi=v_z$ satisfies the boundary value problem
\begin{equation}\label{vz-eq}
-\Delta\psi=\omega_{zz},
~~\hbox{$\psi=0$ on $\partial D_a$.}
\end{equation}
Hence, case (i) of Lemma \ref{lem-laplace} applies to the function
$\psi$.  Moreover, the elliptic estimates for $\omega$, as in
(\ref{h2-est}) imply that
$\|\omega_{zz}\|_{L^2(D_a)}\le\|\Delta\omega\|_{L^2(D_a)}\le
C\|\nabla T\|_{L^2(D_a)}(1+\|\nabla T\|_{L^2(D_a)})$.
Hence, the same proof as in the derivation of
the bound (\ref{omega-infty-bd}) applies to $\psi$ and we obtain that
\[
\|v_z\|_{L^\infty(D_a)}\le C\|\nabla T\|_{L^2(D_a)}(1+\|\nabla T\|_{L^2(D_a)}).
\]
This, together with (\ref{omega-infty-bd}) implies that
\[
\|w_x\|_{L^\infty(D_a)}\le C\|\nabla T\|_{L^2(D_a)}(1+\|\nabla T\|_{L^2(D_a)}).
\]

The other pair of derivatives, $v_x$ and $w_z$, do not satisfy
a homogeneous boundary condition on the lines $x=\pm a$. Therefore, one
cannot apply the standard elliptic estimates up to the boundary to
the function $\eta=w_z=-v_x$ (the second equality follows from the
incompressibility of the flow).  In order to circumvent this
difficulty, we extend
the function $w$ to a larger cylinder
$D_{a+1}=[-a-1,a+1]\times[0,\lambda]$ by setting $w(-a-x,z)=w(-a+x,z)$
and $w(a+x,z)=w(a-x,z)$ for $0\le x\le 1$. The resulting function is
of a class $C^2(D_{a+1})$ since $w(x,z)$ satisfies the Neumann
boundary condition at $x=\pm a$. This also extends the function
$\eta=w_z$ to the larger cylinder. Moreover, $\eta$ satisfies the
Neumann boundary condition along the horizontal lines $z=0,\lambda$:
\[
\eta_z=w_{zz}=-v_{zx}=0\hbox{ on $z=0,\lambda$,}
\]
and
\[
\Delta\eta=\omega_{xz},
\]
with the function $\omega$ extended to the larger cylinder by the same
reflection.  Hence, the interior elliptic estimates up to the boundary
for solutions of the Neumann problem imply that
\[
\|\eta\|_{H^2(Q)}\le\|\Delta\eta\|_{L^2(D_a)}+\|\eta\|_{L^2(D_a)}
\]
for any rectangle $Q=[x_0,x_0+1]\times[0,\lambda]$ that is strictly
contained inside the larger cylinder $D_{a+1}$. Therefore, the Sobolev
embedding theorem together with the above estimates imply that
\begin{equation}\label{eta-infty-0}
\|\eta\|_{L^\infty(D_a)}\le C[\|\Delta\eta\|_{L^2(D_{a})}+
\|\eta\|_{L^2(D_{a})}]=
C[\|\omega_{xz}\|_{L^2(D_a)} +\|\eta\|_{L^2(D_a)}].
\end{equation}
However, as the function $\omega$ satisfies the Dirichlet boundary
conditions in $D_a$, we can apply the estimate (\ref{h2-est}) to the
function $\omega$ up to boundary, to obtain
\[
\|\omega\|_{H^2(D_a)}\le C\left[\|\Delta\omega\|_{L^2(D_a)}+
\|\omega\|_{L^2(D_a)}\right].
\]
We now use the vorticity equation (\ref{vort-eq}) to bound
$\|\Delta\omega\|_{L^2(D_a)}$ and the estimate (\ref{omega-T})
to estimate $\|\omega\|_{L^2(D_a)}$, and conclude that
\begin{equation}\label{omegaxz-T}
\|\omega_{xz}\|_{L^2(D_a)}\le C\|\nabla T\|_{L^2(D_a)}(1+\|\nabla T\|_{L^2(D_a)}).
\end{equation}
Furthermore, the estimate (\ref{nablaw-T}) for $\|\nabla
w\|_{L^2(D_a)}$ implies that
\begin{equation}\label{eta-T}
\|\eta\|_{L^2(D_a)}=\|w_z\|_{L^2(D_a)}\le\|\nabla w\|_{L^2(D_a)}\le
C\|\nabla T\|_{L^2(D_a)}.
\end{equation}
We infer from the bounds (\ref{eta-infty-0}), (\ref{omegaxz-T}) and
(\ref{eta-T}) that
\[
\|\eta\|_{L^\infty(D_a)}\le
C\|\nabla T\|_{L^2(D_a)}(1+\|\nabla T\|_{L^2(D_a)}).
\]
This proves the uniform bound on $w_z$ and hence the proof of Lemma
\ref{lem-uTbd} is complete. $\Box$


Let us now proceed to estimate $\|\nabla T\|_{L^2(D_a)}$ in terms of
$\|v\|_{L^\infty(D_a)}$, a bound in the direction opposite to that in
Lemma \ref{lem-uTbd}. Most importantly, we will bound the square
$\|\nabla T\|_2^2$ in terms of a linear function of $\|v\|_\infty$. As
we are unable to obtain such bound by the standard ellitpic estimates,
we have to proceed with an explicit calculation. As a preliminary step
we show the following.
\begin{lemma}\label{lem-Tu-bd}
Let $(c,T,\vu)$ satisfy (\ref{1.1.tau-2})-(\ref{bc-lat-a.tau-2}) with
the normalization (\ref{T1/2.tau}). Then there exists a constant $C>0$
and a constant $a_0>0$ so that we have for all $a>a_0$ and
$0\le\tau\le 1$
\begin{equation}\label{Tu-bd}
\int_{D_a}|\nabla T|^2dxdz+
\int_0^\lambda T_x(a,z)dz\le C\left[1+\|v\|_\infty\right].
\end{equation}
\end{lemma}
{\bf Proof.} Recall that the function $T$ satisfies
\begin{equation}\label{T-eq-2}
-cT_x+\tau\vu\cdot\nabla T=\Delta T+\tau f(T)
\end{equation}
with the boundary conditions
\begin{equation}\label{bc-T-trav-fin-2}
\hbox{$T(-a,z)=1$,  $T(a,z)=0$,
$\disp\pdr{T}{z}=0$ at $z=0,\lambda$.}
\end{equation}
We multiply (\ref{T-eq-2}) by $(1-T)$ and use the boundary conditions
and incompressibility of the flow $\vu$ to obtain
\begin{equation}\label{c-a-eq}
\frac{c\lambda}{2}=\int_0^\lambda T_x(a,z)dz+\int_{D_a}|\nabla T|^2dxdz
+\tau\int (1-T)f(T)dxdz.
\end{equation}
%
Hence Lemma \ref{lem-cu}, and the fact that $(1-T)f(T)\ge 0$ imply that
\begin{equation}\label{ctau-nabla}
\int_{D_a}|\nabla T|^2dxdz+
\int_0^\lambda T_x(a,z)dz\le\frac{c\lambda}{2}\le
C\left[1+\|v\|_\infty\right]
\end{equation}
and Lemma \ref{lem-Tu-bd} is proved. $\Box$

In order to close the bounds (\ref{ctau-bds}), (\ref{uT-bd}) and
(\ref{Tu-bd}) we need to bound the integral of $T_x$ in
(\ref{Tu-bd}). This is done in the next Lemma.
\begin{lemma}\label{lem-Tx}
Let $(c,T,\vu)$ satisfy
(\ref{1.1.tau-2})-(\ref{bc-lat-a.tau-2}) with the normalization
(\ref{T1/2.tau}). There exists a constant $C>0$ and a constant $a_0$
so that we have for all $a\ge a_0$ and $0\le\tau\le 1$
\begin{equation}\label{Tx-bd}
0\le -\int_0^\lambda T_x(a,z)dz\le
C\left[1+\|\nabla T\|_2\right].
\end{equation}
\end{lemma}
{\bf Proof.} In order to find a bound for $\disp\int_0^\lambda T_x(x=\pm
a,z)dz$ we introduce
\[
I(x)=\frac{1}{\lambda}\int_0^\lambda T(x,z)dz
\]
and integrate equation (\ref{1.1.tau-2}) for $T$ in $z$. Using
the boundary conditions we obtain
\begin{equation}\label{eq-G}
-I_{xx}=G(x),~~I(-a)=1,~~I(a)=0,~~
G(x)=\frac{\tau}{\lambda}\int f(T(x,z))dz
-\int(\tau\vu\cdot\nabla T-c T_x)\frac{dz}{\lambda}.
\end{equation}
This equation can be solved explicitly:
\[
I(x)=-\int_{-a}^x(x-s)G(s)ds+Ax+B
\]
with constants
\[
A=-\frac{1}{2a}+\frac{1}{2a}\int_{-a}^a(a-s)G(s)ds,~~
B=\frac{1}{2}+\frac{1}{2}\int_{-a}^a(a-s)G(s)ds
\]
that are determined from the boundary conditions.
Thus, we have
\[
I_x(-a)=A,~~I_x(a)=A-\int_{-a}^aG(s)ds.
\]
Using the expression for the function $G(x)$ in
(\ref{eq-G}), we now infer that
\begin{eqnarray*}
&&0\le -I_x(a)=\frac{1}{2a}+\frac{1}{2a}\int_{-a}^a(a+s)G(s)ds=
\frac{1}{2a}+\frac{\tau}{2a}\int_{-a}^a\int_0^\lambda(a+x)f(T(x,z))
\frac{dzdx}{\lambda}\\
&&-\frac{\tau}{2a}\int_{-a}^a\int_0^\lambda(a+x)\vu\cdot\nabla
T(x,z)\frac{dzdx}{\lambda}
+\frac{c}{2a}\int_{-a}^a\int_0^\lambda(a+x)T_x
\frac{dzdx}{\lambda}.
\end{eqnarray*}
Integrating by parts, using the boundary conditions and
incompressibility of $\vu$ we obtain
\begin{eqnarray*}
&&0\le -I_x(a)=
\frac{1}{2a}+\frac{\tau}{2}\int_{-a}^a\int_0^\lambda f(T(x,z))
\frac{dzdx}{\lambda}+
\frac{\tau}{2}\int_{-a}^a\int_0^\lambda\frac{x}{a}f(T(x,z))
\frac{dzdx}{\lambda}\\
&&+\frac{\tau}{2a}\int_{-a}^a\int_0^\lambda
v(x,z)T(x,z)\frac{dzdx}{\lambda}
-\frac{c}{2a}\int_{-a}^a\int_0^\lambda
T\frac{dzdx}{\lambda}.
\end{eqnarray*}
However, the normalization condition (\ref{T1/2.tau}) implies that
$f(T(x,z))=0$ for $x\ge 0$ since there is no reaction to the right of
$x=0$. Therefore, we can drop the third term above. This is one of the
crucial points in the proof of the current lemma. Hence, we conclude that
\begin{eqnarray}
&&0\le -I_x(a)
\le\frac{1}{2a}+\frac{\tau}{2}\int_{-a}^a\int_0^\lambda f(T(x,z))
\frac{dzdx}{\lambda}
+\frac{\tau}{2a}\int_{-a}^a\int_0^\lambda
v(x,z)T(x,z)\frac{dzdx}{\lambda}\nonumber\\
&&-\frac{c}{2a}\int_{-a}^a\int_0^\lambda T \frac{dzdx}{\lambda}
\le \frac{1}{2a}+\frac{\tau}{2}\int_{-a}^a\int_0^\lambda f(T(x,z))
\frac{dzdx}{\lambda}+
\tau\|v\|_\infty+|c|.\label{Ixbd}
\end{eqnarray}
We used the fact that $0\le T\le 1$ to bound the last term above.
Next, we look at $I_x(-a)$:
\begin{eqnarray*}
&&0\le -I_x(-a)=\frac{1}{2a}-\frac{1}{2a}\int_{-a}^a(a-s)G(s)ds=
\frac{1}{2a}-\frac{\tau}{2a}\int_{-a}^a\int_0^\lambda (a-x)f(T(x,z))
\frac{dzdx}{\lambda}\\
&&+\frac{\tau}{2a}\int_{-a}^a\int_0^1(a-x)\vu\cdot\nabla
T(x,z)dzdx -\frac{c}{2a}\int_{-a}^a\int_0^\lambda
(a-x)T_x \frac{dzdx}{\lambda}.
\end{eqnarray*}
We can drop the second term above, as $(a-x)f(T)\ge 0$, so that,
after integration by parts, we get
\begin{eqnarray}
&&0\le -I_x(-a)\le\frac{1}{2a}
+\frac{\tau}{2a}\int_{-a}^a\int_0^\lambda v(x,z)T(x,z)
\frac{dzdx}{\lambda}
-\frac{c}{2a}\int_{-a}^a\int_0^\lambda T
\frac{dzdx}{\lambda}+|c|\nonumber\\
&&\le \frac{1}{2a}+\tau\|v\|_\infty+|c|.\label{Ix(-a)bd}
\end{eqnarray}
Let us now put together (\ref{Ixbd}), (\ref{Ix(-a)bd}) and (\ref{c-a-eq})
We observe that, with $\disp F=\tau\int
f(T)\frac{dxdz}{\lambda}$, we have the following three inequalities:
\begin{eqnarray*}
&&c=I_x(a)-I_x(-a)+F\\
&&0\le -I_x(a)\le \frac{1}{2a}+\frac{F}{2}
+\tau\|v\|_\infty+|c|\\
&&0\le -I_x(-a)\le \frac{1}{2a}+\tau\|v\|_\infty+|c|.
\end{eqnarray*}
This implies that
\begin{equation}\label{F-u}
F\le \frac{1}{a}+2\tau\|v\|_\infty+4|c|
\end{equation}
and thus
\[
-I_x(a)\le 3|c|+\frac{1}{a}+2\tau\|v\|_\infty\hbox{ for $a\ge a_0$}
\]
Lemma \ref{lem-cu} implies then that
\[
-I_x(a)\le C\left[1+\tau\|v\|_\infty+\frac{1}{a}\right].
\]
Finally, Lemma \ref{lem-uTbd} implies that
\[
-I_x(a)\le C\left[1+\|\nabla T\|_2+
\frac{1}{a}\right] \hbox{ for $a\ge a_0$}.
\]
Thus, Lemma \ref{lem-Tx} is proved. $\Box$

The previous lemmas imply uniform bounds that we summarize as follows.
\begin{corollary}\label{cor-bds}
Let $(c,T,\vu)$ satisfy
(\ref{1.1.tau-2})-(\ref{bc-lat-a.tau-2}) with the normalization
(\ref{T1/2.tau}). There exists a constant $C>0$ and $a_0>0$ so that we have
for all $a\ge a_0$ and $0\le\tau\le 1$
\begin{eqnarray}\label{bds-a-summ}
\|\vu\|_{L^\infty(D_a)}+\|\nabla\vu\|_{L^\infty(D_a)}+
\|\omega\|_{L^2(D_a)}+|c|+\|\nabla T\|_{L^2(D_a)}
+\tau\int_{D_a}f(T)dxdz\le C.
\end{eqnarray}
\end{corollary}
In particular, as a consequence we also have
\begin{eqnarray}\label{bds-a-summ-2}
\|\nabla\vu\|_{L^2(D_a)}+
\|\omega\|_{L^\infty(D_a)}\le C.
\end{eqnarray}
{\bf Proof.} Lemmas \ref{lem-Tx} and \ref{lem-Tu-bd} imply that
\[
\int|\nabla T|^2dxdz\le C\left[1+\tau\|v\|_\infty+\|\nabla T\|_{L^2(D_a)}\right].
\]
Then Lemma \ref{lem-uTbd} implies that
\[
\|\nabla T\|_{L^2(D_a)}^2\le
C(1+\|\nabla T\|_{L^2(D_a)})
\]
and thus the estimate on $\|\nabla T\|_{L^2(D_a)}$ in
(\ref{bds-a-summ}) holds. Then Lemma
\ref{lem-uTbd} implies the bounds on $\|\vu\|_{L^\infty(D_a)}$,
$\|\nabla\vu\|_{L^\infty(D_a)}$ and $\|\omega\|_{L^2(D_a)}$.  The
bound on $|c|$ in (\ref{bds-a-summ}) now follows from Lemma
\ref{lem-cu}. Finally, the estimate on the total reaction rate follows
from the above bounds and (\ref{F-u}). One can elliminate the factor
$\tau$ in front of the total reaction rate in (\ref{bds-a-summ}):
actually, one can show that it remains bounded as $\tau\to 0$.
However, unlike the other estimates in (\ref{bds-a-summ}), we will use
the bound on the total reaction rate only at $\tau=1$. $\Box$

It remains to prove the uniform H\"older $C^{1,\alpha}$-estimates for $T(x,z)$,
$\omega$ and $\vu$ in order to finish the proof of Proposition
\ref{lem-1/2}.
\begin{lemma}\label{lem-holder}
There exist two constants $C>0$ and $a_0>0$ so that the following bound
holds for all $a\ge a_0$:
\begin{equation}\label{gradT-infty}
\|\omega\|_{C^{1,\alpha}(D_a)}+\|\vu\|_{C^{1,\alpha}(D_a)}
+\|T\|_{C^{1,\alpha}(D_a)}\le C
\end{equation}
provided that $0\le\alpha<1$.
\end{lemma}
{\bf Proof.} The bound for $T$ follows from the standard elliptic
local regularity estimates up to the boundary \cite{GT}, the
$C^1$-bound on the flow $\vu$ and the uniform bound on the speed $c$
in Corollary \ref{cor-bds}. The H\"older estimate for $\omega$ follows
then from the vorticity equation (\ref{vort-eq}) with the Dirichlet
boundary conditions, the above mentioned $C^{1,\alpha}$-bound on $T$,
the same uniform estimates in Corollary
\ref{cor-bds} and the same results of \cite{GT}.
Finally, the H\"older bounds on $\vu$ follow
from the Poisson equations (\ref{u-poisson}) and (\ref{vert-eq}) on
the horizontal and vertical flow components, respectively, and the
H\"older estimates for $\omega$. $\Box$

This completes the proof of the a priori bounds in Proposition
\ref{lem-1/2}. We now turn to the proof of the existence part of this
proposition.

{\bf Step 2. The degree argument.} The a priori bounds proved in the
first step of the proof allow us to use the Leray-Schauder topological
degree argument to establish existence of solutions to the problem
({\ref{1.1.tau-2})-(\ref{bc-lat-a.tau-2}) with the normalization
(\ref{T1/2.tau}) in the bounded domain $D_a$. This method of
construction of traveling wave solutions goes back to
\cite{BNS}. We introduce a map
\[
{\cal K_\tau}:(c,\omega,T)\to (\theta^\tau,\Omega^\tau,Z^\tau)
\]
as the solution operator of the linear system
\begin{eqnarray}\label{1.1.tau}
-cZ_x^\tau+\tau\vu\cdot\nabla Z^\tau=\Delta Z^\tau+\tau f(T)\\
-c\Omega_x^\tau+\vu\cdot\nabla\Omega^\tau-\sigma\Delta\Omega^\tau=
\tau\rho[e_2 T_x-e_1T_z]\nonumber
\end{eqnarray}
in $D_a$ with the no stress boundary conditions
\begin{equation}\label{bc-lat.tau}
\pdr{Z^\tau}{z}=0,~~\tilde w^\tau=\Omega^\tau=0~~\hbox{at $z=0,\lambda$}
\end{equation}
and
\begin{equation}\label{bc-lat-a.tau}
Z^\tau(-a,y)=1,~~Z^\tau(a,y)=0,~~
\tilde u^\tau=\Omega^\tau=0~~\hbox{at $x=\pm a$.}
\end{equation}
Here the unknown flow $\tilde\vu^\tau=(\tilde u^\tau ,\tilde w^\tau)$ and
the given flow $\vu$ are the
incompressible flows corresponding to
the vorticities $\Omega^\tau$ and $\omega$, respectively, and satisfying
the no-stress boundary conditions.
The number $\theta^\tau$ is defined by
\[
\theta^\tau=\theta_0-\max_{x\ge 0} T(x,z)+c.
\]
The operator ${\cal K}_\tau$ is a mapping of the Banach space
$X=\Rm\times C^{1,\alpha}(D_a)\times C^{1,\alpha}(D_a)$, equipped with the norm
$\|(c,\omega,T)\|_X=\max(|c|,\|\omega\|_{C^{1,\alpha}(D_a)},\|T\|_{C^{1,\alpha}(D_a)})$,
onto itself.
A solution $\bq^\tau=(c^\tau,\omega^\tau,T^\tau)$ of
({\ref{1.1.tau-2})-(\ref{bc-lat-a.tau-2}) is a fixed point of ${\cal
K}_\tau$ and satisfies ${\cal K}_\tau\bq^\tau=\bq^\tau$, and vice
versa: a fixed point of ${\cal K}_\tau$ provides a solution to
({\ref{1.1.tau-2})-(\ref{bc-lat-a.tau-2}). Hence, in order to show that
(\ref{1.1.tau-2})-(\ref{bc-lat-a.tau-2}) has a traveling front solution
it suffices to show
that the kernel of the operator ${\cal F}_\tau=\hbox{Id}-{\cal K}_\tau$ is not
trivial. The standard elliptic regularity results in \cite{GT} imply
that the operator ${\cal K}_\tau$ is compact and depends continuously
on the parameter $\tau\in[0,1]$. Thus the Leray-Schauder
topological degree theory can be applied. Let us introduce a ball
$B_M=\{\|(c,\omega,T)\|_X\le M\}$. Then Lemma
\ref{lem-holder} and Lemma \ref{lem-cu} show that the operator ${\cal F}_\tau$
does not vanish on the boundary $\partial B_M$ with $M$ sufficiently
large for any $\tau\in[0,1]$. It remains only to show that the degree
$\hbox{deg}({\cal F}_1,B_M,0)$ in $\bar B_M$ is not zero. However, the
homotopy invariance property of the degree implies that
$\hbox{deg}({\cal F}_\tau,B_M,0)=\hbox{deg}({\cal F}_0,B_M,0)$ for all
$\tau\in[0,1]$. Moreover, the degree at $\tau=0$ can be computed
explicitly as the operator ${\cal F}_0$ is given by
\[
{\cal F}_0(c,\omega,T)=(\max_{x\ge 0} T(x,y)-\theta_0,\omega,T-T_0^c).
\]
Here the function $T_0(x)$ solves
\[
\frac{d^2T_0^c}{dx^2}+c\frac{dT_0^c}{dx}=0,~~~T_0^c(-a)=1,~~T_0^c(a)=0
\]
and is given by
\[
T_0^c(x)=\frac{e^{-cx}-e^{-ca}}{e^{ca}-e^{-ca}}.
\]
The mapping ${\cal F}_0$ is homotopic to
\[
\Phi(c,\omega,T)=(\max_{x\ge 0} T_0^c(x,y)-\theta_0,\omega,T-T_0^c)
\]
that in turn is homotopic to
\[
\tilde\Phi(c,\omega,T)=
(T_0^c(0)-\theta_0,\omega,T-T_0^{c_*^0}),
\]
where $c_*^0$ is the unique number so that $T_0^{c_*}(0)=\theta_0$.
The degree of the mapping $\tilde\Phi$ is the product of the degrees
of each component.  The last two have degree equal to one, and the
first to $-1$, as the function $T_0^c(0)$ is decreasing in $c$. Thus
$\hbox{deg}{\cal F}_0=-1$ and hence $\hbox{deg}{\cal F}_1=-1$ so that
the kernel of $\hbox{Id}-{\cal K}_1$ is not empty. This finishes the
proof of Proposition \ref{lem-1/2}. $\Box$

\begin{remark} Observe that the
$C^{1,\alpha}$-regularity of $T$, $\vu$ and
$\omega$ can be bootstrapped to $C^{2,\alpha}$-regularity: we have
\begin{equation}\label{gradT-infty-c2}
\|\omega\|_{C^{2,\alpha}(D_{a})}+\|\vu\|_{C^{2,\alpha}(D_{a})}
+\|T\|_{C^{2,\alpha}(D_{a})}\le C
\end{equation}
provided that $0\le\alpha<1$.
\end{remark}

\section{Identification of the limit}\label{sec:limit}

In order to finish the proof of Theorem \ref{thm-main} we consider
the solutions $(c^a,T^a,\vu^a)$ constructed in Proposition
\ref{lem-1/2} and pass to the limit $a\to +\infty$.
The a priori estimates in the same proposition imply that we can
choose a subsequence $a_n\to\infty$ so that $T_n(x,z)=T_{a_n}(x,z)$
converges uniformly on compact sets to a function $T(x,z)$, while the
flow $\vu_n(x,z)=\vu_{a_n}(x,z)$ converges to a flow $\vu(x,z)=(v,w)$ and
the front speeds also converge: $c_n=c_{a_n}\to c$. The vorticity
functions $\omega_n(x,z)=\omega_{a_n}(x,z)$ converge to the limit
$\omega=w_x-v_z$. The limits satisfy the uniform bounds
\begin{equation}\label{T-bd}
|c|+\|\vu\|_\infty+\|\omega\|_\infty+
\int|\nabla T|^2dxdz+\int f(T)dxdz+\int|\nabla\vu|^2dxdz
\le C
\end{equation}
that follow from Corollary \ref{cor-bds} and the H\"older estimates
(\ref{gradT-infty-holder}) and (\ref{gradT-infty-c2}).
The regularity estimates on $(T^a,\vu^a)$
imply that the limit functions $T$ and $\vu$ satisfy the Boussinesq
system
\begin{eqnarray}\label{T-eq}
&&-cT_x+\vu\cdot\nabla T=\Delta T+f(T)\nonumber\\
&&-c\omega_x+\vu\cdot\nabla\omega-\sigma\Delta\omega=
\rho(\hat\ve\cdot\nabla^\perp T)\nonumber\\
&&\omega=w_x-v_z.
\end{eqnarray}
Moreover, the boundary conditions on the lateral boundaries hold
for $T$ and $\vu$:
\begin{equation}\label{3-bc}
\pdr{T}{z}=0,~~w=\omega=0\hbox{ on $z=0,\lambda$.}
\end{equation}
The normalization condition
\begin{equation}\label{normal-limit}
~~\max_{x\ge 0} T(x,z)=\theta_0
\end{equation}
is also satisfied.

Therefore, to finish the proof of Theorem \ref{thm-main}, it remains
only to show that (i) $T$ converges to a constant $\theta_-$ as $x\to
-\infty$ and $T\to 0$ as $x\to +\infty$, (ii) $\vu\to 0$ as
$x\to+\infty$, and (iii) $\theta_-=1$ if the reaction rate satisfies
$f(T)\le (T-\theta_0)_+^2/\lambda^2$.  First, we note that the uniform
$L^2$-bound on $\nabla T$ in (\ref{T-bd}) implies that $T$ converges
to two constants $\theta_-$ and $\theta_+$ as $x\to\pm\infty$,
possibly passing to a subsequence $x_n\to\pm\infty$. The elliptic
regularity results imply that actually $T$ converges to these
constants as $x\to\pm\infty$. Moreover, the bound for the total
reaction rate $\disp\int f(T)dxdz$ in (\ref{T-bd}) implies that
$f(\theta_-)=f(\theta_+)=0$. Furthermore, integrating (\ref{T-eq}) we
obtain
\begin{equation}\label{c-int}
c(\theta_--\theta_+)=\int f(T)\frac{dxdz}{\lambda}.
\end{equation}

In order to identify the limits $\theta_\pm$ we will make use of the
following lemmas that provide some additional information on solutions
on a finite domain before the passage to the limit. The first result
describes the behavior near the right end $x=a_n$.
\begin{lemma}\label{lem:T_x(a)}
There exists a sequence $a_n\to\infty$ so that
\begin{equation}\label{Txn0}
\disp\left|\pdr{T_n(a_n,z)}{x}\right|\to 0
\end{equation}
as $n\to\infty$, uniformly in $z$. Moreover, we
have $\lim_{n\to\infty} T_n(a_n-x_0,z)=0$ for all $x_0\in\Rm$.
\end{lemma}
{\bf Proof.}  We introduce a shifted solution
$\Phi_n(x,z)=T_n(x+a_n,z)$ , $\bv_n=\vu_n(x+a_n,z)$ defined in the
domain $-2a_n\le x\le 0$.  The functions $\Phi_n$ and $\bv_n$ satisfy
the same a priori bounds (\ref{T-bd}) as $T_n$ and $\vu_n$ and hence
they converge as $n\to\infty$ to some limits $\Phi$ and $\bv$ that
satisfy
\begin{equation}\label{Phi-eq}
-c\Phi_x+\bv\cdot\nabla \Phi=\Delta \Phi,~~\Phi(0,z)=0, ~~x\le 0
\end{equation}
as $f(\Phi_n)=0$ for $x>-a_n$ and thus in the limit $f(\Phi)=0$. The
function $\Phi$ satisfies the Neumann boundary conditions at
$z=0,\lambda$. The uniform upper bound on $\|\nabla\Phi_n\|_2$
together with the elliptic regularity results imply that $\Phi$ has to
converge to a constant $\Phi_-$ as $x\to-\infty$ along a
subsequence. We note that, as $0\le\Phi(x,z)\le\theta_0$, the constant
$\Phi_-$ satisfies the same bounds:
\[
0\le\Phi_-\le\theta_0.
\]
Integrating (\ref{Phi-eq}) we obtain
\begin{equation}\label{c-Phi-}
c\lambda\Phi_-=\int \Phi_x(0,z)dz\le 0.
\end{equation}
Hence, either $\Phi_-=0$ or $c\le 0$. In the former case $\Phi\equiv
0$ and hence $\Phi_x(0,z)=0$ for all $z$. That implies that both
$T_x^n(a_n,z)\to 0$ as $n\to\infty$ and
$T_x^n(a_n-x_0,z)=\Phi_n(-x_0,z)\to 0$ as $n\to\infty$, as claimed in
Lemma \ref{lem:T_x(a)}. It remains to rule out the second case, $c\le
0$. This is done in the next lemma that provides a crucial lower bound
on the speed $c_n$. In particular it shows that $c>0$ -- this will
conclude the proof of Lemma \ref{lem:T_x(a)}.
\begin{lemma} \label{lem-c-lower}
The front speed is positive, $c>0$.
\end{lemma}
{\bf Proof.} Integrating the temperature equation in (\ref{T-eq}) for
$T_n$, we obtain
\begin{equation}\label{cn-eq}
c_n\lambda=\int_0^\lambda \pdr{T_n}{x}(a_n,z)dz-
\int_0^\lambda \pdr{T_n}{x}(-a_n,z)dz+
\int f(T_n)dxdz \ge\int_0^\lambda \pdr{T_n}{x}(a_n,z)dz+\int f(T_n)dxdz.
\end{equation}
Observe also that we have a uniform bound
\begin{equation}\label{prod-bd}
\left(\int f(T_n)dxdz\right)\left(\int|\nabla T_n|^2dxdz\right)\ge C
\end{equation}
that follows from the fact that $T_n(0,z)\le \theta_0$
and $T_n(-a_n,z)=1$.
The proof of (\ref{prod-bd}) is as in \cite{CKOR}:
there exists  $z_0\in(0,\lambda)$ such that both
\[
\int\limits_{-a_n}^0 |\nabla T_n(x,z_0)|^2dx \leq
3\int\limits_{D_n}|\nabla T_n|^2\frac{dxdz}{\lambda},~~
D_n=[-a_n,a_n]_x\times[0,\lambda]_z,
\]
and
\[
\int\limits_{-a_n}^0 f(T_n(x,z_0))dx\le 3\int\limits_{D_n}f(T_n(x,z))
\frac{dxdz}{\lambda}.
\]
Let $x_1$ be the left-most point so that
$\disp T_n(x_1,z_0)=1-\frac{1-\theta_0}{4}$:
\[
x_1=\inf\left\{x\in(-a_n,0):~T_n(x,z_0)=1-\frac{1-\theta_0}{4}\right\}
\]
and $x_2>x_1$ be the left-most point so that
$\disp T_n(x_2,z_0)=\theta_0+\frac{1-\theta_0}{4}$:
\[
x_2=\inf\left\{x\in(-a_n,0):~T_n(x,z_0)=
\theta_0+\frac{1-\theta_0}{4}\right\}.
\]
Existence of $x_1$ and $x_2$ is guaranteed by the fact that
$T_n(-a_n,z)=1$ and $T_n(0,z)\le\theta_0$ for all $z\in[0,\lambda]$.
Then the reaction rate $f(T_n(x,z_0))>C$ for $x_1\le x\le x_2$ so that
\[
C|x_1-x_2|\le\int_{x_1}^{x_2} f(T_n(x,z_0))dx\le
3\int\limits_{D_n}f(T_n(x,z))\frac{dxdz}{\lambda}
\]
and
\[
\frac{(1-\theta_0)^2}{4|x_1-x_2|}\le
\int_{x_1}^{x_2}|\nabla T_n(x,z_0)|^2dx\le
 3\int\limits_{D_n}|\nabla T_n|^2\frac{dxdz}{\lambda}.
\]
Multiplying these two inequalities, we arrive at (\ref{prod-bd}).

The estimate (\ref{prod-bd}) and the uniform upper bound on $\|\nabla
T_n\|_2$ in Corollary \ref{cor-bds} imply that
\begin{equation}\label{lower-fT}
\int_{D_n} f(T_n)dxdz\ge C.
\end{equation}
Then, passing to the limit in (\ref{cn-eq}), and using (\ref{c-Phi-})
we obtain
\[
c\lambda(1-\Phi_-)\ge C>0,
\]
as
\[
\int \pdr{T_n}{x}(a_n,z)dz\to\int\Phi_x(0,z)dz,
\]
with the function $\Phi$ as in the proof of Lemma
\ref{lem:T_x(a)}. Now, we recall that $\Phi_-\le\theta_0<1$ and thus
the front speed $c>0$.  This finishes the proof of Lemma
\ref{lem-c-lower} and hence also that of Lemma \ref{lem:T_x(a)}. $\Box$

Lemma \ref{lem-c-lower} and (\ref{c-int}) imply that
$\theta_-\ge\theta_+$. However, if $\theta_-=\theta_+$ we have
$f(T)=0$ everywhere and hence (\ref{T-eq}) is a linear equation. The
maximum principle implies that $T\equiv\hbox{const}$ in this case. The
last condition in (\ref{T-bd}) implies that this constant has to be
equal to $\theta_0$. Hence, either $\theta_->\theta_+$ or
$T\equiv\theta_0$.

Let us now rule out the special case that $\theta_-=\theta_+=\theta_0$.
\begin{lemma}\label{lem-theta+-}
The left and right limits $\theta_-$ and $\theta_+$ satisfy
$\theta_->\theta_+$.
\end{lemma}
{\bf Proof.} We have already shown
that $\theta_-\ge\theta_+$ and, moreover, if $\theta_-=\theta_+$ then
\begin{equation}\label{3-ridic}
\theta_-=\theta_+=\theta_0.
\end{equation}
Hence, it suffices to show that the latter
is impossible. Let us assume that (\ref{3-ridic}) holds.
As we have explained above, then
\begin{equation}\label{3-ridic2}
\hbox{$T_n\to\theta_0$, and
$\partial T_n/\partial x\to 0$ uniformly on compact sets}.
\end{equation}
Then, integrating the equation
\[
-c_n\pdr{T_n}{x}+\vu_n\cdot\nabla T_n=\Delta T_n+f(T_n)
\]
between $x=0$ and $x=a_n$ we obtain, as $f(T_n)=0$ in this region,
\begin{equation}\label{3-theta0?}
c_n\int_0^\lambda T_n(0,z)dz-\int_0^\lambda v_n(0,z)T_n(0,z)dz=
\int_0^\lambda \pdr{T_n}{x}(a_n,z)dz-
\int_0^\lambda \pdr{T_n}{x}(0,z)dz.
\end{equation}
We now pass to the limit $n\to\infty$ in (\ref{3-theta0?}). The first
term on the left converges to $c\theta_0\lambda$, as we have assumed
that $T$ converges uniformly to $\theta_0$ on compact intervals. The
second term on the left converges to
\[
\int_0^\lambda v(0,z)\theta_0dz=0,
\]
as incompressibility of the flow $\vu_n$ and the boundary conditions
at $x=\pm a$ imply that
\[
\int v_n(0,z)dz=0.
\]
The limit (\ref{Txn0}) in Lemma \ref{lem:T_x(a)} implies that the
first term on the right side of (\ref{3-theta0?}) converges to
zero. Finally, the last term on the right side of (\ref{3-theta0?})
converges to zero because of (\ref{3-ridic2}).  Therefore, we obtain
\[
c\lambda\theta_0=0.
\]
However, this implies that $c=0$ which contradicts Lemma
\ref{lem-c-lower}. Hence, the case $\theta_-=\theta_+=\theta_0$ is
ruled out and thus $\theta_->\theta_+$. $\Box$

We continue the analysis of the behavior of the solution at the right
end of the domain.
\begin{lemma}\label{lem-grad}
The gradient $\nabla T_n$ converges to zero ``as $x\to +\infty$''
uniformly in $n$, that is, for every $\eps>0$ there exists
$N\in{\mathbb N}$ and $R$ so that $|\nabla T_n(x,z)|<\eps$ for all
$n\ge N$ and all $R\le x\le a_n$.
\end{lemma}
{\bf Proof.} Let us assume that this is not the case. Then there
exists $\eps_0>0$ and a sequence $b_n\to +\infty$ so that $|\nabla
T_n(b_n,z_n)|\ge\eps_0$ for some $z_n\in[0,\lambda]$. Note that Lemma
\ref{lem:T_x(a)} implies that
\begin{equation}\label{anbn}
|a_n-b_n|\to\infty.
\end{equation}
Let us define the shifted solution $\Psi_n(x,z)=T_n(x-b_n,z)$,
$\bv_n(x,z)=\vu_n(x-b_n,z)$ on the domain $x\in[-a_n-b_n,a_n-b_n]$.
Then $\Psi_n$ and $\bv_n$ satisfy the same uniform bounds as $T_n$ and
$\vu_n$ and thus they converge to a pair of functions $\Psi$, $\bv$
uniformly on compact intervals, together with their derivatives. The
functions $\Psi$ and $\bv$ are defined on the whole real line because
of (\ref{anbn}).  Moreover, the function $\Psi$ has left and right
limits $\Psi_\pm$ as $x\to\pm\infty$. Hence, the same argument as in
the proof of Lemma \ref{lem:T_x(a)} shows that $\Psi$ must be equal a
constant, as it has left and right limits and satisfies
\[
-c\Psi_x+\bv\cdot\nabla\Psi=\Delta\Psi.
\]
However, this contradicts the fact that
$\max_z|\nabla\Psi(0,z)|\ge\eps_0$. $\Box$

The decay of the gradient of $T_n$ implies that the flow ahead of the
front goes to zero for large $x$, uniformly in $n$.
\begin{lemma}\label{lem-u} The
flow $\vu(x,z)$ converges to zero
on the right uniformly in $n$, that is, for any $\eps>0$ there
exists $R>0$ and $N\in{\mathbb N}$ so that $|\vu(x,z)|<\eps$ for
all $R<x\le a_n$.
\end{lemma}
{\bf Proof.} We choose $N$ and $R$ so that $|\nabla T_n|<\eps$ for all
$n\ge N$ and $x\in[R,a_n]$. Then, we decompose
$T_n=T_n^{in}+T_n^{out}$ with $\hbox{supp}~T_n^{in}\subset\{x\le
R+1\}$ and $\hbox{supp}~ T_n^{out}\subset\{x\ge R\}$. We also require
that both $T_n^{in}$ and $T_n^{out}$ satisfy the same uniform gradient
bounds as $T_n$. Moreover, we have $|\nabla T_n^{out}|<\eps$. We also
split $\omega_n=\omega_n^{in}+\omega_n^{out}$ and
$\vu_n=\vu_n^{in}+\vu_n^{out}$ accordingly:
\[
-c_n\omega_n^{in}+\vu_n\cdot\nabla\omega_n^{in}-\sigma\Delta\omega_n^{in}=
\rho\left(e_2\pdr{T_n^{in}}{x}-e_1\pdr{T_n^{in}}{z}\right),
~~\omega_n^{in}=0\hbox{ on $\partial D_{a_n}$},
\]
and similarly for $\omega_n^{out}$.

We now bound $|\vu_n^{in}|$ and
$|\vu_n^{out}|$ separately for sufficiently large $x$. First, we
look at $\vu_n^{in}=(v_n^{in},w_n^{in})$. The function $\omega_n^{in}$
satisfies a homogeneous equation
\begin{equation}\label{3-hom-eq}
-c_n\pdr{\omega_n^{in}}{x}+\vu_n\cdot\nabla\omega_n^{in}-
\sigma\Delta\omega_n^{in}=0
\end{equation}
in the rectangle $D_{R+2,a_n}=\{R+2\le x\le a_n,~0\le z\le\lambda\}$,
as $T_n^{in}$ vanishes in $D_{R+2,a_n}$. The function $\omega_n^{in}$
satisfies a uniform $C^{2,\alpha}$-bound -- this is shown in the same
way as the $C^{2,\alpha}$-bound for the full vorticity function
$\omega$ in (\ref{gradT-infty-c2}). This in turn implies that the
function $\psi(z)=\omega_n^{in}(R+2,z)$ is uniformly bounded in
$C^2[0,\lambda]$. Let $g(x)$ be a smooth monotonic and positive
cut-off function so that
\begin{equation}\label{3-fn-g}
\hbox{$g(x)=1$ for $R+2\le x\le R+3$ and $g(x)=0$ for
$x\ge R+4$.}
\end{equation}
Then the function $\omega_n^{in}$ can be decomposed as
\[
\omega_{n}^{in}(x,z)=\psi(z)g(x)+\zeta_n.
\]
The function $\zeta_n$ satisfies
\begin{equation}\label{3-zeta-eq}
-c_n\pdr{\zeta_n}{x}+\vu_n\cdot\nabla\zeta_n
-\sigma\Delta\zeta_n=f_n  \hbox{ in $D_{R+2,a_n}$,
$\zeta_n=0$ on $\partial D_{R+2,a_n}$.}
\end{equation}
The right side $f_n$ is given by
\[
f_n:=\sigma\psi''(z)g(x)+\sigma\psi(z)g''(x)-c_n\psi(z)g'(x)
-v_n\psi(z)g'(x)-w_n\psi'(z)g(x).
\]
It is
supported in $R+2\le x\le R+4$ and is uniformly
bounded since $\|\omega_n^{in}\|_{C^{2,\alpha}(D_a)}\le C$. Let us
choose $\alpha>0$ sufficiently small, then the function
$\xi_n(x,z)=\zeta_n(x,z)e^{\alpha x}$ satisfies
\begin{equation}\label{3-xi-eq}
-c_n\pdr{\xi_n}{x}+\alpha c\xi_n+\vu_n\cdot\nabla\xi_n-\alpha v_n\xi_n
-\sigma\Delta\xi_n+2\sigma\alpha\pdr{\xi_n}{x}-\sigma\alpha^2\xi_n=g_n
 \hbox{ in $D_{R+2,a_n}$,
$\xi_n=0$ on $\partial D_{R+2,a_n}$}
\end{equation}
with $g_n=f_n(x)e^{\alpha x}$.
Multiplying (\ref{3-xi-eq}) by $\xi_n$ and integrating by parts,
using the boundary conditions, we obtain
\begin{equation}\label{xi-bd}
\sigma\int_{D_{R+2,a_n}}|\nabla
\xi_n|^2dxdz+\left(c\alpha-\alpha\|v\|_\infty-\sigma\alpha^2\right)
\int_{D_{R+2,a_n}}|\xi_n|^2dxdz\le\|g_n\|_2\|\xi_n\|_2.
\end{equation}
However, as the function $\xi_n$ vanishes at $z=0,\lambda$,
the Poincar\'e inequality implies that
\[
\int_{D_{R+2,a_n}}|\nabla \xi_n|^2dxdz\ge\frac{\pi^2}{\lambda^2}
\int_{D_{R+2,a_n}}|\xi_n|^2dxdz.
\]
Hence, the following upper bound holds
\[
\int_{D_{R+2,a_n}}|\xi_n|^2dxdz\le\|g_n\|_2^2\le C,
\]
provided that $\alpha$ is sufficiently small, since $\|v\|_\infty\le C$.
Using (\ref{xi-bd}) once again we conclude that
\[
\int_{D_{R+2,a_n}}|\nabla\xi_n|^2dxdz \le C.
\]
Therefore, the function
$\zeta_n$ satisfies
\[
\int_{D_{R+2,a_n}}\left[|\nabla\zeta_n|^2+|\zeta_n|^2\right]
 e^{2\alpha x}dxdz\le C.
\]
This, in turn implies the same bound for the function
$\omega_{n}^{in}$:
\begin{equation}\label{3-omegain-l2}
\int_{R+2}^{a_n}\int_0^\lambda|\omega_n^{in}|^2 e^{2\alpha x}dxdz
+\int_{R+2}^{a_n}\int_0^\lambda|\nabla\omega_n^{in}|^2
e^{2\alpha x}dxdz\le C.
\end{equation}
It follows that the $L^2$-norm of $\omega_n^{in}$ decays uniformly
in $n$:
\begin{equation}\label{3-omegain-l2-unif}
\int_{r_0}^{a_n}\int_0^\lambda|\omega_n^{in}|^2dxdz\le
e^{-\alpha r_0}\int_{r_0}^{a_n}\int_0^\lambda|\omega_n^{in}|^2
e^{2\alpha x}dxdz
\le C e^{-\alpha r_0}.
\end{equation}
for $r_0>R+5$, and the same bound holds for $\nabla\omega_n^{in}$:
\begin{equation}\label{3-nabla-omegain-l2-unif}
\int_{r_0}^{a_n}\int_0^\lambda|\nabla\omega_n^{in}|^2dxdz\le
e^{-\alpha r_0}\int_{r_0}^{a_n}\int_0^\lambda|\nabla\omega_n^{in}|^2
e^{2\alpha x}dxdz
\le C e^{-\alpha r_0}.
\end{equation}

As the function $\omega_{n}^{in}$ satisfies the homogeneous equation
(\ref{3-hom-eq}) for $x\ge R+1$ with a bounded flow $\vu$, the
standard local elliptic estimates now imply that
\begin{equation}\label{3-omegain-unif}
|\omega_n^{in}(x,z)|\le Ce^{-\alpha x} \hbox{ for $x\ge R+5$}.
\end{equation}
The $W^{2,p}$ elliptic estimates imply then the uniform decay of
the gradient of $\omega_n^{in}$:
\begin{equation}\label{3-omegain-grad}
|\nabla\omega_n^{in}(x,z)|\le Ce^{-\alpha x} \hbox{ for $x\ge R+5$}.
\end{equation}

Now we can bound the flow $\vu_n^{in}=(v_n^{in},w_n^{in})$ itself.
First, we look at the horizontal component $v_n^{in}$. It
satisfies the following Poisson equation in $D_{R+2,a_n}$:
\[
-\Delta v_n^{in}=\pdr{\omega_n^{in}}{z}\hbox{ in $D_{R+2,a_n}$, }
\pdr{v_n^{in}}{z}=0\hbox{ on $z=0,\lambda$, }
v_n^{in}=0\hbox{ on $x=a_n$.}
\]
Moreover, the $C^{2,\alpha}$-regularity of $\vu_n$ implies that
the boundary value $\phi(z)=v_n^{in}(R+2,z)$ is bounded in
$C^2[0,\lambda]$. Therefore, as we did with $\omega_{n}^{in}$,
we represent $v_{n}^{in}(x,z)=\phi(z)g(x)+\bar v_n^{in}(x,z)$
with the cut-off function $g(x)$ as in (\ref{3-fn-g}). The
function $\bar v_n^{in}$ satisfies
\[
-\Delta \bar v_n^{in}=\bar f_n:=-\phi_{zz}(z)g(x)-\phi(z)g''(x)
+\pdr{\omega_n^{in}}{z}\hbox{ in $D_{R+2,a_n}$, }
\]
with an exponentially decaying function $\bar f_n$, as follows
from (\ref{3-omegain-grad}).
The boundary conditions are
\[
\pdr{\bar v_n^{in}}{z}=0\hbox{ on $z=0,\lambda$, }
\bar v_n^{in}=0\hbox{ on $x=R+2,a_n$.}
\]
The same argument as we used to obtain (\ref{3-omegain-l2})
implies that
\begin{equation}\label{3-vin-l2}
\int_{R+2}^{a_n}\int_0^\lambda|v_n^{in}|^2 e^{2\beta x}dxdz
\le C\int|\bar f_n(x,z)|^2e^{2\beta x}dxdz\le C
\end{equation}
with a sufficiently small $0<\beta<\alpha$.
Therefore, in the same vein as we have obtained (\ref{3-omegain-unif})
and (\ref{3-omegain-grad}), we conclude that
\begin{equation}\label{3-vin-unif}
|v_n^{in}(x,z)|\le Ce^{-\alpha x} \hbox{ for $x\ge R+5$}
\end{equation}
and
\begin{equation}\label{3-vin-grad}
|\nabla v_n^{in}(x,z)|\le Ce^{-\alpha x} \hbox{ for $x\ge R+5$}.
\end{equation}
The uniform bound on $w_n^{in}$ now follows, as it satisfies the
Dirichlet boundary condition $w_n^{in}(x,0)=w_n^{in}(x,\lambda)=0$ and
the derivative $\displaystyle\pdr{w_n^{in}}{z}= -\pdr{v_n^{in}}{x}$ is
exponentially decaying (\ref{3-vin-grad}). We infer that
\begin{equation}\label{3-win-unif}
|w_n^{in}(x,z)|\le Ce^{-\alpha x} \hbox{ for $x\ge R+5$}.
\end{equation}

Now we bound $\vu_n^{out}$. The corresponding vorticity satisfies
\[
-c_n\pdr{\omega_n^{out}}{x}+\vu_n\cdot\nabla\omega_n^{out}-
\sigma\Delta\omega_n^{out}=\rho(\hat\ve\cdot\nabla^\perp
T_n^{out})\hbox{ in $D_{a_n}$,
$\omega=0$ on $D_{a_n}$.}
\]
However,
\begin{equation}\label{3-Tout-eps}
|\nabla T_n^{out}|\le \eps
\end{equation}
by construction, hence the maximum principle implies that
\begin{equation}\label{3-omegaout-bd}
|\omega_n^{out}(x,z)|\le \eps \rho q(z)\le C\eps.
\end{equation}
Here the non-negative function $q(z)$ satisfies the boundary value
problem
\[
-\sigma q''(z)+w_n(z)q'(z)=1,~~q(0)=q(\lambda)=0.
\]
We infer from the standard local elliptic estimates up to
the boundary, (\ref{3-omegaout-bd})  and (\ref{3-Tout-eps}) that
\begin{equation}\label{3-nablaomega-out}
|\nabla\omega_n^{out}(x,z)|\le  C\eps\hbox{ in $D_{a_n}$}
\end{equation}
as well. The vertical flow component satisfies
\[
\Delta w_n^{out}=\pdr{\omega_n^{out}}{x}\hbox{ in $D_{a_n}$,
$w_n^{out}(x,0)=w_n^{out}(x,\lambda)=0,
\displaystyle\pdr{w_n^{out}(\pm a_n,z)}{x}=0.$}
\]
Therefore, the maximum principle implies once again that
\[
|w_n^{out}(x,z)|\le \frac{C\rho\eps}{2}z(\lambda-z)\le C\eps\hbox{ in
$D_{a_n}$.}
\]
Hence, the same local elliptic regularity results allow us to
conclude that
\[
|\nabla w_n^{out}(x,z)| \le C\eps\hbox{ in $D_{a_n}$.}
\]
In order to bound the horizontal flow component $v_n^{out}$ and
conclude the proof of Lemma \ref{lem-u} we observe that
$\displaystyle\pdr{v_n^{out}}{z}=\pdr{w_n^{out}}{z}-\omega_n^{out}$ so
that $\displaystyle\left|\pdr{v_n^{out}}{z}\right|\le C\eps$ in
$D_{a_n}$. However, $v_n^{out}$ also satisfies the mean-zero condition
\[
\int_0^\lambda v_n^{out}(x,z)=0\hbox{ for all $-a_n\le x\le a_n$.}
\]
Hence, we have $|v_n^{out}(x,z)|\le C\eps$  in $D_{a_n}$, and the
proof of Lemma \ref{lem-u} is now complete. $\Box$


The next lemma implies that the right limit $\theta_+=0$.
\begin{lemma}\label{lem-right}
The right limit $\theta_+=0$.
\end{lemma}
{\bf Proof.} Let us choose $R$ independent of $n$ so that
$c_{a_n}>\sup_{x>R}|v_n(x,z)|$ for all $n$.  Lemma \ref{lem-c-lower}
implies that the speeds $c_n$ are uniformly bounded below by a
positive constant, thus it follows from Lemma \ref{lem-u} that we can
find such $R>0$. Then the function $\phi(x)=Ae^{-\alpha x}$, with a
sufficiently small $\alpha>0$, satisfies
\[
-c_n\phi_x+\vu_n\cdot\nabla\phi\ge \Delta\phi.
\]
An argument as in the proof of Lemma \ref{lem-cu} shows that if $A$ is
chosen so that $Ae^{-\alpha R}>1$ then $T_n(x,z)\le A e^{-\alpha x}$
on the domain $x\in[R,a_n]$. Therefore, the limit $T(x,z)$ obeys the
same bound, which in turn implies that $\theta_+=0$. $\Box$

Finally, we show that under the additional assumption
(\ref{f-technical}) the left limit $\theta_-=1$. This is the only
place in the proof where assumption (\ref{f-technical}) is used.
\begin{lemma}\label{lem-left-one}
Let us assume that $f(T)\le(T-\theta_0)_+^2/\lambda^2 $.  Then
the left limit $\theta_-=1$.
\end{lemma}
{\bf Proof.} We note that we have for each $x\in\Rm$
\[
\int |\nabla T(x,z)|^2dz\ge \frac{(M(x)-m(x))^2}{\lambda}
\]
with $M(x)=\max_z T(x,z)$ and $m(x)=\min_zT(x,z)$. It follows from
the maximum principle  that the function $m(x)$ is non-increasing.
Let us assume that $\theta_-\le\theta_0$,
then  monotonicity of $m(x)$ implies that $m(x)<\theta_0$ for all
$x\in\Rm$. Then we have
\[
\int |\nabla T(x,z)|^2dxdz\ge \int(M(x)-m(x))^2\frac{dx}{\lambda}\ge
\int (T(x,z)-\theta_0)_+^2\frac{dxdz}{\lambda^2}.
\]
We also observe that
\[
c\theta_-=\int f(T)\frac{dxdz}{\lambda},~~
\frac{c\theta_-^2}{2}+\int|\nabla T|^2\frac{dxdz}{\lambda}=
\int Tf(T)\frac{dxdz}{\lambda}
\]
so that
\[
\int|\nabla T|^2dxdz=\int \left(T-\frac{\theta_-}{2}\right)f(T)dxdz.
\]
Hence we obtain using (\ref{f-technical})
\[
\int\left(T-\frac{\theta_-}{2}\right)(T-\theta_0)_+^2
\frac{dxdz}{\lambda^2}\ge
\int \left(T-\frac{\theta_-}{2}\right)f(T)dxdz\ge
\int (T(x,z)-\theta_0)_+^2\frac{dxdz}{\lambda^2}.
\]
However,  the left side is smaller than the right side unless
$T\equiv\theta_0$, the case that we have already ruled
out.  $\Box$

This finishes the proof of Theorem \ref{thm-main}.

\section{Bounds for the initial value problem}\label{sec:bds}

We consider in this section the solutions of the Cauchy problem with
general front-like initial data and obtain the uniform bounds on the
bulk burning rate and other average quantities stated in Theorems
\ref{intro-thm-bds} and Theorem \ref{intro-thm-narrow}. We
prove the first result, and the proof of the second result  is
presented in Section \ref{sec-narrow}.

\subsection{Bounds in an arbitrary strip}

We prove in this section Theorem
\ref{intro-thm-bds}. Let $T(t,x,z)$, $u(t,x,z)$ be the solution of the Cauchy problem
\begin{eqnarray}\label{cauchy}
T_t+\vu\cdot\nabla T=\Delta T+f(T)\\
\vu_t+\vu\cdot\nabla\vu-\sigma\Delta\vu+\nabla p=
\rho T\hat\ve\nonumber\\
\nabla\cdot \vu=0,\nonumber\\
\end{eqnarray}
with initial data $T_0(x,z)$, $\vu_0(\vx)$.  We assume that there
exists $R>0$ so that $T_0(x,z)=0$ for $x>R$ and $T_0(x,z)=1$ for
$x<-R$, and that the initial vorticity is bounded in $L^2$:
\[
\int|\omega_0(x,z)|^2dxdz<+\infty.
\]
The assumptions on the initial temperature $T_0$ can be relaxed --  it
simply has to approach one and zero at the two ends of the
domain sufficiently fast.

We recall that the bulk burning rate $\bar V(t)$, the Nusselt number
$\bar N(t)$ and the average horizontal flow $\bar U(t)$ are defined by
\begin{eqnarray}\label{4-def}
&&\bar V(t)=\frac{1}{t}\int_0^t\int V(s)ds,~~
V(t)=\int f(T)\frac{dxdz}{\lambda},\\
&&\bar N(t)=\frac{1}{t}\int_0^tN(s)ds,~~N(t)=\int |\nabla
T|^2\frac{dxdz}{\lambda},\\
&&\bar U(t)=\frac{1}{t}\int_0^t\|v(s)\|_\infty ds.
\end{eqnarray}
The laminar front speed $c_0$ is defined as the unique $c$ so that
equation
\[
-c\Phi'=\Phi''+f(\Phi),~~\Phi(-\infty)=1,~\Phi(+\infty)=0
\]
has a solution $0<\Phi<1$. We recall the statement of Theorem
\ref{intro-thm-bds}.
\begin{thm}\label{thm-bds}
There exists a constant $C>0$ so that
under the above assumptions on the initial data $T_0,\vu_0$, the following bounds hold
\begin{eqnarray}\label{4-bds}
&&c_0-C[\rho+\rho^2]+o(1)\le\bar V(t)\le c_0+C[\rho+\rho^2]+o(1)\\
&&\bar N(t)\le \left[{C\rho}+
\sqrt{\frac{c_0}{2}+C^2\rho^2}\right]^2+o(1)\nonumber\\
&&\bar U(t)\le C\rho[1+\rho]+o(1)\nonumber
\end{eqnarray}
as $t\to +\infty$.
\end{thm}
This theorem  shows that the coupled problem (\ref{cauchy})
is in a sense a regular perturbation of the single
reaction-diffusion equation with $\rho=0$.
The lower bound in (\ref{4-bds}) is of interest only for small
$\rho$ when the left side is positive.

{\bf Proof.} First, we prove the following bounds on $\bar N(t)$ and
$\bar V(t)$ in terms of $\bar U(t)$.
\begin{lemma}\label{4-lem-barN-bd}
There exists a constant $C_0$ that depends only on the initial data
$T_0$ so that
\begin{equation}\label{4-VU-bd-0}
\bar N(t)\le \frac{1}{2}\bar V(t)+\bar U(t)+
C_0\left[\frac{1}{{t}}+\frac{1}{\sqrt{t}}\right]
\end{equation}
and
\begin{equation}\label{kansas-0}
\bar V(t)\le c_0+\bar U(t)+C_0\left[\frac{1}{t}+\frac{1}{\sqrt{t}}\right].
\end{equation}
\end{lemma}
{\bf Proof.}
Define $g(T)=T(1-T)$
and its integral
\[
R(t)=\int g(T)\frac{dxdz}{\lambda}.
\]
The idea of using a concave function $g(T)$ in a related context is
due to B. Winn \cite{Winn}. We observe that
\begin{equation}\label{bliams}
\frac{dR}{dt}=\int g'(T)\Delta T\frac{dxdz}{\lambda}+
\int g'(T)f(T)\frac{dxdz}{\lambda}\ge
-\int g''(T)|\nabla T|^2\frac{dxdz}{\lambda}-V(t)
\end{equation}
with the burning rate $V(t)$ defined in (\ref{4-def}).
Thus
\[
\frac{dR}{dt}+V(t)\ge 2\int|\nabla T|^2\frac{dxdz}{\lambda}=2N(t),
\]
which after  averaging in time becomes
\begin{equation}\label{bums}
\frac{R(t)}{t}+\bar V(t)\ge 2\bar N(t).
\end{equation}

In order to obtain an upper bound for the potentially small term
$R(t)/t$ in (\ref{bums}) we construct sub- and super-solutions for
$T(t,x,z)$. This construction follows \cite{X2}. We look for a
sub-solution for $T$ of the form
\[
\psi_l(t,x,z)=\Phi_0(x-c_0t+x_1+\xi_1(t))-q_1(t,x,z).
\]
Here $\Phi_0$ is the traveling wave in the absence of convection, at
$\rho=0$, normalized so that $\Phi_0(0)=\theta_0$. It is the
unique solution of
\[
-c_0\Phi_0'=\Phi_0''+f(\Phi_0), ~~\Phi_0(0)=\theta_0,~
\Phi_0(-\infty)=1,~\Phi_0(+\infty)=0.
\]
The functions $\xi_1(t)$ and $q_1(t,x,z)$ are to be chosen.
In order for $\psi_l$ to be a sub-solution we need
\[
G[\psi_l]=\pdr{\psi_l}{t}+
\vu\cdot\nabla\psi_l-\Delta\psi_l-f(\psi_l)\le 0.
\]
We have
\[
G[\psi_l]=\dot{\xi}_1\Phi_0'+u\Phi_0'-
\pdr{q_1}{t}-\vu\cdot\nabla q_1+\Delta
q_1+f(\Phi_0)-f(\Phi_0-q_1).
\]
With an appropriate choice of $x_1$, that is, by shifting $\Phi_0$
sufficiently to the left we can ensure that
$T_0(x,z)\ge\Phi_0(x)-q_{10}(x)$ with $0\le q_{10}(x)\le
(1-\theta_0)/2$ and $q_{10}(x)\in L^1(\Rm)$.  Then we choose
$q_1(t,x,z)$ to be the solution of
\begin{equation}\label{4-adv-diff}
\pdr{q_1}{t}+\vu\cdot\nabla q_1=\Delta q_1,~~q_1(0,x,z)=q_{10}(x),
\pdr{q_1}{z}=0\hbox{ at $z=0,\lambda$}.
\end{equation}
The following lemma first proved in \cite{FKR}
provides a uniform $L^1-L^\infty$ decay estimate
for $q_1$ that is independent of the advection term.
\begin{lemma}\label{lem-FKR} There exists a constant $C>0$ that is
independent of the (incompressible) flow $\vu$ so that
\begin{equation}\label{4-to-zero}
\|q_1(t)\|_\infty\le\frac{C}{\lambda\sqrt{t}}\|q_{10}\|_{L^1(D)}
\end{equation}
for $t\ge 1$.
\end{lemma}
As mentioned above, the main point of the above result is the
independence of the constant in (\ref{4-to-zero}) from the flow
$\vu$. We also note that this $L^1-L^\infty$ estimate behaves in a
one-dimensional way for large times, as one would expect for a
strip. The factor of $\lambda$ in the denominator is compensated by
the fact that the $L^1$-norm is taken over the strip and not only in
$x$. We postpone the proof of Lemma \ref{lem-FKR} till the end of this
section.

We can find $\delta>0$ so that if $\Phi_0\in(1-\delta,1)$ and
$q_1\in(0,(1-\theta_0)/2)$ then $f(\Phi_0)\le f(\Phi_0-\delta)$. Hence
we have in this range of $\Phi_0$:
\begin{equation}\label{4-ineq-N}
G[\psi_l]\le \dot{\xi}_1\Phi_0'+v\Phi_0'.
\end{equation}
Furthermore, if $\Phi_0\in (0,\delta)$ then
$f(\Phi_0)=f(\Phi_0-\delta)=0$ and hence in this range of $\Phi_0$ we
have (\ref{4-ineq-N}) with the equality sign. Finally, if
$\Phi_0\in(\delta,1-\delta)$ then $|f(\Phi_0)-f(\Phi_0-q)|\le
K|q|$ and $\Phi_0'\le -\beta$. Hence $G[\psi_l]\le 0$ everywhere
provided that
\begin{equation}\label{4-ineq-xi}
\dot\xi_1(t)\ge\|v(t)\|_\infty+\frac{K\|q(t)\|_\infty}{\beta}.
\end{equation}
Thus choose
\begin{equation}\label{4-xi}
\xi_1(t)=\bar U(t)t+C\sqrt{t}.
\end{equation}
Therefore we obtain a lower bound for $T$:
\begin{equation}\label{4-T-lbd}
T(t,x,z)\ge \Phi_0(x-c_0t+\bar U(t)t+C\sqrt{t})-q_1(t,x,z).
\end{equation}
In order to obtain an upper bound we set
$\psi_u=\Phi_0(x-c_0t-x_2-\xi_2(t))+q_2(t,x,z)$ and look for
$\xi_2(t)$ and $q_2(t,x,z)$ so that $G[\psi_u]\ge 0$.  The constant
$x_2$ is chosen so that
\[
T_0(x,z)\le \Phi_0(x-x_2)+q_2(0,x,z)
\]
with $q_2(0,x,z)\in L^1(D)$ and $0\le q_2(0,x,z)\le
\theta_0/2$, as with $q_1(0,x,z)$.
The function $q_2(t,x,z)$ is then chosen to satisfy the same
advection-diffusion equation (\ref{4-adv-diff}) similarly to
$q_1$. Hence it obeys the same time decay bounds as $q_1$. With the
above choice of $q_2$ we have
\[
G(\psi_u)=-\dot{\xi}_2\Phi_0'+v\Phi_0'
+f(\Phi_0)-f(\Phi_0+q_2).
\]
Once again, we consider three regions of values for $\Phi_0$.
First, if $1-\delta\le\Phi_0\le 1$ with a sufficiently small $\delta>0$
then $f(\Phi_0)-f(\Phi_0+q_2)\ge 0$, as $q_2\ge 0$. Hence $G[\psi_u]\ge 0$ in this region
provided that $\dot\xi_2\ge 0$. Second, as $q_2\le\theta_0/2$
we have $f(\Phi_0)=f(\Phi_0+q_2)=0$ if $0\le\Phi_0\le\delta$ with
a sufficiently small $\delta>0$. Hence $G[\psi_u]\ge 0$ in that region
under the same condition $\dot\xi_2\ge 0$. Finally, if
$\Phi_0\in(\delta,1-\delta)$ then $\Phi_0'\le -\beta$ with
$\beta>0$ and $|f(\Phi_0)-f(\Phi_0+q_2)|\le K\|q_2\|_\infty$.
That means that $G[\psi_u]\ge 0$ if we choose $\xi_2$ so that
\[
\dot\xi_2\ge \|v(t)\|_\infty+\frac{K\|q_2\|_\infty}{\beta}.
\]
Therefore we choose
\[
\xi_2(t)=\bar U(t)t+C\sqrt{t},
\]
as with $\xi_1(t)$. Therefore we obtain upper and lower bounds
\begin{equation}\label{4-T-lbd-ubd-0}
\Phi_0(x-c_0t+\xi_1(t)+x_1)-q_1(t,x,z)\le
T(t,x,z)\le\Phi_0(x-c_0t-\xi_2(t)-x_2)+q_2(t,x,z)
\end{equation}
that imply in particular that
\begin{equation}\label{4-T-lbd-ubd}
\Phi_0(x-c_0t+\bar U(t)t+C_0[1+\sqrt{t}])-\frac{C_0}{\sqrt{t}}
\le T(t,x,z)\le
\Phi_0(x-c_0t-\bar U(t)t-C_0[1+\sqrt{t}])+\frac{C_0}{\sqrt{t}}
\end{equation}
with a constant $C_0$ determined by the initial conditions.
Hence, using (\ref{4-T-lbd-ubd-0})-(\ref{4-T-lbd-ubd}) and the
$L^1$-bounds $\|q_j(t)\|_{L^1(D)}\le C_0$, $j=1,2$, we obtain
\begin{eqnarray*}
&&R(t)=\int T(1-T)\frac{dxdz}{\lambda}=
\int_{-\infty}^{c_0t-\xi_2(t)-x_2}\int
T(1-T)\frac{dzdx}{\lambda}+
\int_{c_0t-\xi_2(t)-x_2}^{c_0t+\xi_1(t)+x_1}\int
T(1-T)\frac{dzdx}{\lambda}\\
&&+\int_{c_0t+\xi_1(t)+x_1}^{\infty}\int
T(1-T)\frac{dzdx}{\lambda}
\le C_0+\int_{-\infty}^{0}
(1-\Phi_0)dx+(\xi_1(t)+\xi_2(t))+\int_{0}^{\infty}
\Phi_0(x)dx\\
&&\le C_0(1+\sqrt{t})+2t\bar U(t).
\end{eqnarray*}
This together with (\ref{bums}) implies that
\begin{equation}\label{4-VU-bd}
\bar V(t)+2\bar U(t)+C_0\left[\frac{1}{{t}}+\frac{1}{\sqrt{t}}\right]
\ge 2\bar N(t)
\end{equation}
so that (\ref{4-VU-bd-0}) holds.

Moreover, we have
\begin{eqnarray}\label{kansas}
&&\bar V(t)=\frac{1}{t}\int_0^t\left(\int f(T(s,x,z))
\frac{dxdz}{\lambda}\right)ds=
\frac{1}{t}\int_0^t\left(\int T_t(s,x,z)\right)\frac{dxdz}{\lambda}ds\\
&&=
\frac{1}{t}\int
\left[T(t,x,z)-T_0(x,z)\right]\frac{dxdz}{\lambda}\nonumber\\
&&\le \frac{1}{t}
\int\left[\Phi_0(x-c_0t-x_2-\xi_2(t))+q_2(t,x,z)-\Phi_0(x+x_1)+q_1(0,x,z)
\right]
\frac{dxdz}{\lambda}\nonumber\\
&&\le
\frac{C_0}{t}+\frac{1}{t}\int_{-\infty}^{c_0t+\xi_2(t)+x_2}(1-\Phi_0(x))dx+
\frac{1}{t}\int_0^\infty\Phi_0(x)dx\le
c_0+\bar U(t)+C_0\left[\frac{1}{\sqrt{t}}+\frac{1}{t}\right]\nonumber
\end{eqnarray}
as follows from (\ref{4-T-lbd-ubd-0})-(\ref{4-T-lbd-ubd}). This proves
(\ref{kansas-0}) and finishes the proof of Lemma
\ref{4-lem-barN-bd}. $\Box$

On the other hand we have the following upper bound for $\bar U(t)$ in terms
of $\bar N(t)$.
\begin{lemma}\label{lem-4-W-bd}
There exists a constant $C>0$ so that for all $t>0$ the following
inequality holds
\begin{equation}\label{california-0}
\bar U(t)\le C\left[\rho\sqrt{ \bar N(t)}+
\frac{1}{\sqrt{t}}\|\omega_0\|_{L^2}\right].
\end{equation}
\end{lemma}
{\bf Proof.} We multiply the vorticity equation
\[
\pdr{\omega}{t}+\vu\cdot\nabla\omega-\sigma\omega=
\rho(\hat\ve\cdot\nabla^\perp T)
\]
by $\omega$ and integrate:
\begin{equation}\label{4-55}
\frac{1}{2}\frac{d}{dt}\int|\omega(t,x,z)|^2dxdz+
\sigma\int|\nabla\omega(t,x,z)|^2dxdz=\rho\int \omega(t,x,z)
(\hat\ve^\perp\cdot\nabla T)dxdz,
\end{equation}
with $\hat\ve^\perp=(e_2,-e_1)$.

\commentout{
Let us define $\overline
T(t,x,z;\hat\ve^\perp)$ as the average of $T(t,x,z)$ in the direction
$\hat\ve^\perp$ perpendicular to gravity direction $\hat\ve$. That is,
\[
\overline T(t,x,z;\hat\ve^\perp)=\frac{1}{l(\hat\ve^\perp)}
\int_{-\infty}^{\infty} T(t,x+se_2,z-se_1)ds,
\]
where $T$ is set to be zero outside the strip $D$ and $l(\hat\ve^\perp)=
\lambda/e_1$ is the length of the segment that connects the lines $z=0$
and $z=\lambda$ and goes in direction $\hat\ve^\perp$. Note that
$\hat\ve^\perp\cdot\nabla \overline T=0$ so that (\ref{4-55}) can be
re-written as
\begin{eqnarray}\label{4-90}
&&\frac{1}{2}\frac{d}{dt}\int|\omega(t,x,z)|^2dxdz+
\sigma\int|\nabla\omega(t,x,z)|^2dxdz\\
&&~~~~~~~~~=-\rho\int (\hat\ve^\perp\cdot\nabla\omega(t,x,z))
(T(t,x,z)-\overline T(t,x,z;\hat\ve^\perp))dxdz\nonumber\\
&&~~~~~~~~~\le \frac{\sigma}{2}\int|\nabla\omega(t,x,z)|^2dxdz+
C\rho^2\int|T(t,x,z)-\overline T(t,x,z;\hat\ve^\perp)|^2dxdz.\nonumber
\end{eqnarray}
}
The Poincar\'e inequality for $\omega(t,x,z)$ implies then that
\[
\frac{1}{2}\frac{d}{dt}\int|\omega(t,x,z)|^2dxdz+
\frac{\sigma}{2}\int|\nabla\omega(t,x,z)|^2dxdz\le
C\rho^2\int|\nabla T(t,x,z)|^2dxdz.
\]
Integrating this equation in time we conclude that
\begin{equation}\label{4-280}
\frac{1}{t}\int_0^t\int|\nabla\omega(s,x,z)|^2dxdzds
\le C\left[\rho^2\bar N(t)+\frac{1}{t}\|\omega_0\|_{L^2}^2\right].
\end{equation}
However, as in the proof of Lemma \ref{lem-uTbd}, we have
$\|v(t)\|_{L^\infty(D)}\le C\|\nabla\omega(t)\|_{L^2(D)}$. This, together
with (\ref{4-280}) implies (\ref{california-0}). $\Box$


Putting the  bounds (\ref{california-0}) and
(\ref{4-VU-bd-0})-(\ref{kansas-0})
together and using the Cauchy-Schwartz inequality
we arrive at
\[
2\bar N(t)\le
c_0+C_0\left[\frac{1}{t}+\frac{1}{\sqrt{t}}\right]+
\frac{C\rho }{t}\int_0^t\sqrt{  N(s)}ds\le
c_0+C_0\left[\frac{1}{t}+\frac{1}{\sqrt{t}}\right]
+C \rho\sqrt{ \bar N(t)}.
\]
Hence we obtain an upper bound
\begin{equation}\label{4-101}
\bar N(t)\le \left[{C {\rho}}+\sqrt{\frac{c_0}{2}+C^2\rho^2}\right]^2+o(1).
\end{equation}
This, together with (\ref{california-0}) implies that
\[
\bar U(t)\le C\rho\left[1+{\rho}\right]+o(1).
\]
It follows then from
(\ref{kansas}) that
\[
\bar V(t)\le c_0+
C\rho\left[1+\rho\right]+o(1).
\]
The lower bound on $\bar V(t)$ in (\ref{4-bds}) is proved
similarly.
This finishes the proof of Theorem
\ref{thm-bds}. It remains only to prove Lemma \ref{lem-FKR}.
$\Box$

{\bf Proof of Lemma \ref{lem-FKR}.} We will show that there
exists a universal constant $C>0$ so that the solution of
\begin{eqnarray}\label{v-eq}
&&\pdr{\psi}{t}+\vu\cdot\nabla\psi=\sigma\Delta\psi\\
&&\psi(0,x,z)=\psi_0(x,z)\ge 0,\nonumber
\end{eqnarray}
with the Neumann boundary conditions at $z=0$ and $z=\lambda$, and $\vu$
sufficiently regular, satisfies
\begin{equation}\label{l1-linfty}
\|\psi(t)\|_{L^\infty(D)}\le
{C}n^2(t)\|\psi_0\|_{L^1(D)}.
\end{equation}
Here $n(t)$ is the unique solution of
\begin{equation}\label{n-eq}
\frac{n^4(t)}{1+n^3(t)\lambda^3}=\frac{C}{\sigma \lambda^2t}.
\end{equation}

We multiply (\ref{v-eq}) by $\psi$ and integrate over the domain $D$ to obtain
\begin{equation}\label{2.9}
\frac{1}{2}\frac{d}{dt}\|\psi\|_2^2=-\sigma\|\nabla\psi\|_2^2.
\end{equation}
We now prove the following version of the Nash
inequality \cite{Nash} for a strip of width $\lambda$ in two dimensions:
\begin{equation}\label{nash-strip}
\|\nabla\psi\|_2^2\ge
C\frac{\lambda^2\|\psi\|_2^6}{\|\psi\|_1^4+\lambda^3\|\psi\|_1\|\psi\|_2^3}.
\end{equation}
The proof of (\ref{nash-strip}) is similar to that of the usual Nash
inequality. We represent $\psi$ in terms of its Fourier series-integral:
\[
\psi(x,z)=\sum_{n\ge 0}\int_{\mathbb R} e^{ikx}
\cos\left(\frac{\pi nz}{\lambda}\right)
\hat \psi_n(k)\frac{dk}{2\pi},
\]
where
\[
\hat \psi_n(k)=\frac{2}{\lambda}
\int_D e^{-ikx}\cos\left(\frac{\pi nz}{\lambda}\right)\psi(x,z)dxdz.
\]
Therefore we have
\begin{equation}\label{4-sofa}
|\hat \psi_n(k)|\le \displaystyle\frac{2}{\lambda}\|\psi\|_{L^1}.
\end{equation}
The Plancherel formula becomes
\begin{eqnarray*}
&&\int_D|\psi(x,z)|^2dxdz=\sum_{n,m\ge 0}
\int_{{\mathbb R}^2} e^{ikx-ipx}
\cos\left(\frac{\pi n z}{\lambda}\right)\cos\left(\frac{\pi m z}{\lambda}\right)
\hat \psi_n(k)\hat \psi_m^*(p)\frac{dkdpdxdz}{(2\pi)^2}\\
&&~~~~~~~~~~~~~~~~~~~~~~~=
\frac{\lambda}{2}\sum_{n\ge 0}\int|\hat \psi_n(k)|^2\frac{dk}{2\pi}
\end{eqnarray*}
and similarly
\begin{eqnarray*}
&&\int_D|\nabla\psi(x,z)|^2dxdz=
\frac{\lambda}{2}\sum_{n\ge 0}\int_{\mathbb R}\left(k^2+\frac{\pi^2n^2}{\lambda^2}\right)
|\hat \psi_n(k)|^2\frac{dk}{2\pi}.
\end{eqnarray*}
Let $\rho>0$ be a positive number to be chosen later. Then using the above Plancherel
formula we write
\[
\|\psi\|_2^2=I+II,
\]
with the first term that is bounded using (\ref{4-sofa})
\[
I=\frac{\lambda}{2}\sum_{0\le n\le \rho \lambda}
\int_{|k|\le \rho}|\hat \psi_n(k)|^2\frac{dk}{2\pi}\le
\frac{C\lambda\rho([\lambda\rho]+1)}{\lambda^2}\|\psi\|_1^2\le
\frac{C\rho(\lambda\rho+1)}{\lambda}\|\psi\|_1^2.
\]
The rest is  bounded by
\[
II\le \frac{C\lambda}{\rho^2}
\sum_{n\ge 0}\int_{k\in{\mathbb R}}
\left(k^2+\frac{4\pi^2n^2}{\lambda^2}\right)
|\hat\psi_n(k)|^2{dk}\le \frac{C}{\rho^2}\|\nabla\psi\|_2^2.
\]
Therefore we have for all $\rho>0$:
\[
\|\psi\|_2^2\le  \frac{C\rho(\lambda\rho+1)}{\lambda}\|\psi\|_1^2+
\frac{C}{\rho^2}\|\nabla\psi\|_2^2.
\]
We choose $\rho$ so that
\[
\rho^3=\frac{\lambda\|\nabla\psi\|_2^2}{\|\psi\|_1^2}
\]
and obtain
\begin{eqnarray*}
&&\|\psi\|_2^2\le \frac{C\|\nabla\psi\|_2^{2/3}}{\lambda^{2/3}\|\psi\|_1^{2/3}}
\left(\frac{\lambda^{4/3}\|\nabla\psi\|_2^{2/3}}{\|\psi\|_1^{2/3}}+1\right)
\|\psi\|_1^2+\frac{C\|\nabla\psi\|_2^{2}\|\psi\|_1^{4/3}}{\lambda^{2/3}
\|\nabla\psi\|_2^{4/3}}\\
&&~~~~~=
\frac{C}{\lambda^{2/3}}\|\psi\|_1^{4/3}\|\nabla\psi\|_2^{2/3}+
C\lambda^{2/3}\|\nabla\psi\|_{2}^{4/3}\|\psi\|_1^{2/3}.
\end{eqnarray*}
This is a quadratic inequality $ax^2+bx-c\ge 0$ with
$x=\|\nabla\psi\|_2^{2/3}$, $a=C\lambda^{2/3}\|\psi\|_1^{2/3}$,
$b=\displaystyle\frac{C}{\lambda^{2/3}}\|\psi\|_1^{4/3}$,
and $c=\|\psi\|_2^2$ and hence
\[
x\ge\frac{-b+\sqrt{b^2+4ac}}{2a}=\frac{2c}{b+\sqrt{b^2+4ac}}\ge
\frac{c}{\sqrt{b^2+4ac}}.
\]
This implies that
\[
\|\nabla\psi\|_2^{2/3}\ge
C{\|\psi\|_2^2}\left({{\frac{\|\psi\|_1^{8/3}}{\lambda^{4/3}}+
\lambda^{2/3}\|\psi\|_1^{2/3}\|\psi\|_2^2}}\right)^{-1/2}
\]
and therefore
\begin{eqnarray*}
&&\|\nabla\psi\|_2^2\ge C{\|\psi\|_2^6}\left({{\frac{4\|\psi\|_1^{8/3}}{\lambda^{4/3}}+
4\lambda^{2/3}\|\psi\|_1^{2/3}\|\psi\|_2^2}}\right)^{-3/2}\ge
C{\|\psi\|_2^6}\left(\frac{\|\psi\|_1^{4}}{\lambda^{2}}+
\lambda\|\psi\|_1\|\psi\|_2^3\right)^{-1} \\
&&~~~~~~~~
\ge \frac{C\lambda^2\|\psi\|_2^6}{\|\psi\|_1^{4}+\lambda^3\|\psi\|_1\|\psi\|_2^3}.
\end{eqnarray*}
Hence (\ref{nash-strip}) indeed holds.

We insert (\ref{nash-strip}) into  (\ref{2.9}) and use
the conservation of the $L^1$-norm of $\psi$ (recall that the initial
data is non-negative) obtain
\begin{equation}\label{2.26}
\frac{d\|\psi\|_2}{dt}\le -\frac{C\sigma \lambda^2\|\psi\|_2^5}{\|\psi_0\|_1^4+
\lambda^3\|\psi_0\|_1\|\psi\|_2^3}.
\end{equation}
Integrating (\ref{2.26}) in time we have
\[
C\sigma \lambda^2t\le \frac{\|\psi_0\|_1^4}{\|\psi\|_2^4}+
\frac{\lambda^3\|\psi_0\|_1}{\|\psi\|_2}\le
\frac{1}{z(t)}\left[\lambda^3+\frac{1}{z^3(t)}\right],
\]
where $z(t)=\|\psi(t)\|_2/\|\psi_0\|_1$, and thus
\begin{equation}\label{z-ineq}
\frac{z^4(t)}{1+\lambda^3z^3(t)}\le \frac{1}{C\sigma \lambda^2t}.
\end{equation}
The function on the left side of (\ref{z-ineq}) is monotonically increasing
and hence we have
\begin{equation}\label{l1-l2-n}
\|\psi(t)\|_{2}\le n(t)\|\psi_0\|_1,
\end{equation}
where $n(t)$ is the solution of (\ref{n-eq}).

Let us denote by ${\cal P}_t$ the solution operator for (\ref{v-eq}):
$\psi(t)={\cal P}_t\psi_0$. Then (\ref{l1-l2-n}) implies that
$\|{\cal P}_t\|_{L^1\to L^2}\le n(t)$. The adjoint
operator ${\cal P}_t^*$ is the solution operator for
\begin{eqnarray}\label{v-eq*}
&&\pdr{\tilde\psi}{t}-\vu\cdot\nabla\tilde\psi=\sigma\Delta\tilde\psi\\
&&\tilde\psi(0,x)=\tilde\psi_0(x),~~~x\in D\nonumber
\end{eqnarray}
with the Neumann boundary conditions at $z=0,\lambda$.
Note that the preceding estimates rely only on the anti-symmetry of
the convection operator $\vu\cdot\nabla$. Therefore we have the bound
$\|{\cal P}_t^*\|_{L^1\to L^2}\le n(t)$ and hence
$\|{\cal P}_t\|_{L^2\to L^\infty}\le n(t)$ so that
\[
\|\psi(t)\|_{L^\infty}\le n(t/2)\|\psi(t/2)\|_{L^2}
\le  n^2(t/2)\|\psi_0\|_{L^1}
\]
and thus (\ref{l1-linfty}) indeed holds.

The estimate (\ref{4-to-zero}) follows from the observation
that for large
$t\gg 1$, when $n(t)$ is small, we have the bound
\[
n^2(t)\le\frac{C}{(\sigma t)^{1/2}\lambda}.
\]
This finishes the proof of Lemma \ref{lem-FKR}. $\Box$

\subsection{Bounds on the burning rate in a narrow
domain}\label{sec-narrow}

We recall that no non-planar traveling fronts do exist in the
reactive Boussinesq problem in a narrow vertical strip when gravity
is sufficiently small \cite{TV1,TV2,CKR}. Moreover, solutions with general
front-like initial data become asymptotically planar in the long
time limit \cite{CKR}. We extend now this result to the inclined
cylinders. More precisely, we have the following result (this is a re-statement
of Theorem \ref{intro-thm-narrow}).
\begin{thm}\label{thm-narrow} Let $\hat\ve=(e_1,e_2)$ be the unit
vector in the direction of gravity and let $\rho_j=\rho e_j$,
$j=1,2$. There exist two constants $\lambda_0$ and $\rho_0$ so
that if the domain is sufficiently narrow: $\lambda\le\lambda_0$
and gravity is sufficiently small: $\rho\le\rho_0$ then the
burning rate is bounded by
\begin{equation}\label{5-barv}
\bar V(t)\le c_0+C\rho_2+o(1)\hbox{ as $t\to +\infty$.}
\end{equation}
Moreover, the front is nearly planar in the sense that
\begin{equation}\label{5-Tx}
\bar N_z(t)=\frac{1}{t}\int_0^t\|T_z(s)\|_2^2ds\le C\rho_2^2+o(1)
\hbox{ as $t\to +\infty$.}
\end{equation}
\end{thm}
The key point in Theorem \ref{thm-narrow} is that the bounds in
(\ref{5-barv}) and (\ref{5-Tx}) are independent of the gravity
strength $\rho_1$ in the direction parallel to the cylinder.

{\bf Proof.} Multiplying the vorticity equation by $\omega$ and
integrating by parts we obtain
\begin{equation}\label{5-utro1}
\frac{1}{2}\frac{d}{dt}\int|\omega|^2dxdz+\sigma\int|\nabla\omega|^2dxdz=
\rho_2\int T_x\omega dxdz-\rho_1\int T_z\omega dxdz.
\end{equation}
The Poincar\'e inequality applies to $\omega(x,z)$
with the Poincar\'e constant proportional to $1/\lambda$.  Hence,
if $\lambda<\lambda_0$ and $\rho_j<\rho_0$, (\ref{5-utro1}) implies that
\begin{equation}\label{5-utro2}
\frac{1}{2}\frac{d}{dt}\int|\omega|^2dxdz+\frac{\sigma}{2}\int|\nabla\omega|^2dxdz\le
C\rho_2^2\int |T_x|^2  dxdz+C\rho_1^2\int |T_z|^2dxdz.
\end{equation}
We now differentiate the equation for $T(t,x,z)$ in $z$ to get
\[
\pdr{T_z}{t}+\vu\cdot\nabla T_z+\vu_z\cdot\nabla T=\Delta
T_z+f'(T)T_z.
\]
Multiplying this equation by $T_z$ we obtain
\begin{equation}\label{utro-3}
\frac{1}{2}\frac{d}{dt}\int|T_z|^2dxdz+\int|\nabla T_z|^2dxdz+
\int T_z\vu_z\cdot\nabla T dxdz=\int f'(T)T_z^2\le M\int T_z^2dxdz.
\end{equation}
The last integral on the left side is bounded by
\[
\left|\int T_z\vu_z\cdot\nabla T dxdz\right|=\left|\int T\vu_z\cdot\nabla T_z
dxdz\right|\le 2\int |\vu_z|^2dxdz +\frac{1}{2}|\nabla
T_z|^2 dxdz.
\]
This, together with incompressibility of $\vu$, the Poincar\'e inequality
for $T_z$ and (\ref{utro-3})
imply that
\begin{equation}\label{utro-4}
\frac{1}{2}\frac{d}{dt}\int|T_z|^2dxdz+\frac{1}{4}\int|\nabla
T_z|^2dxdz\le 4\int |\omega|^2 dxdz,
\end{equation}
provided that $\lambda<\lambda_0$. Combining (\ref{5-utro2}) and
(\ref{utro-4}) and using the Poincar\'e inequality for $\omega$ and
$T_z$ once again,
we obtain the following inequalities for
$\Omega(t)=\|\omega(t)\|_2^2$ and $N_z(t)=\|T_z\|_2^2$:
\[
\frac{1}{2}\frac{d\Omega}{dt}+\frac{C}{\lambda^2}\Omega\le
C\rho_2^2N_x(t)+C\rho_1^2N_z(t)
\]
and
\[
\frac{1}{2}\frac{dN_z}{dt}+\frac{C}{\lambda^2}N_z\le
4\Omega.
\]
Hence, the function $Q=N_z+\Omega$ satisfies
\[
\frac{1}{2}\frac{dQ}{dt}+\left[\frac{C}{\lambda^2}-4-C\rho_1^2\right]Q\le
C\rho_2^2N_x(t).
\]
Therefore, we have
\[
Q(t)\le Q_0 e^{-\gamma t}+C\rho_2^2\int_0^t e^{-\gamma(t-s)}N_x(s)ds
\]
with $\gamma>0$ provided that $C/\lambda^2>5$ and $C\rho_1^2<1$. We
conclude that
\begin{eqnarray*}
&&\bar Q(t):=\frac{1}{t}\int_0^tQ(\tau)d\tau\le \frac{Q_0}{t}
\int_0^t e^{-\gamma \tau}d\tau+\frac{C\rho_2}{t}\int_0^t \int_0^\tau
e^{-\gamma(\tau-s)}N_x(s)dsd\tau\\
&&\le \frac{C_0}{t}+\frac{C\rho_2^2}{t}\int_0^t N_x(s) e^{\gamma
s}\int_s^te^{-\gamma \tau}d\tau ds\le
\frac{C_0}{t}+\frac{C\rho_2^2}{\gamma}\bar N_x(t)\le
C\rho_2^2+\frac{C_0}{t}.
\end{eqnarray*}
The last inequality above follows from the bound on $\bar N(t)$ in
Theorem \ref{thm-bds}. Now, the bound (\ref{5-Tx}) in Theorem
\ref{thm-narrow} follows. Then, (\ref{5-utro2}) together with
(\ref{5-Tx}) and the same uniform bound on $\bar N(t)$ in Theorem \ref{thm-bds}
imply that
\[
\frac{1}{t}\int_0^t\|\nabla\omega(s)\|_2^2ds\le
C\rho_2^2+\frac{C_0}{t}.
\]
We recall that $\|v(t)\|_{L^\infty(D)}\le
C\|\nabla\omega\|_{L^2(D)}$ -- this, together with the above,
imply that
\begin{equation}\label{5-baru}
\bar U(t)\le C\rho_2+\frac{C_0}{\sqrt{t}}.
\end{equation}
Finally, using (\ref{kansas}) and (\ref{5-baru}) we obtain
(\ref{5-barv}). $\Box$


\begin{thebibliography}{99}

\bibitem{ACVV1} M.~Abel, A.~Celani, D.~Vergni and A.~Vulpiani, \it
Front propagation in laminar flows, \rm Physical Review E, {\bf
64} 6307 (2001).

\bibitem{ACVV2} M.~Abel, M.~Cencini, D.~Vergni and A.~Vulpiani,
\it Front speed enhancement in cellular flows, \rm Chaos {\bf 12},
p. 481.



\bibitem{ABP} B.~Audoly,  H.~Berestycki and Y.~Pomeau, {
    R\`eaction diffusion en \`ecoulement stationnaire rapide},
  C.R.Acad. Sci., Ser. IIB, {\bf 328}, 255-262.

\bibitem{Volpert-private} M. Belk, B. Kazmierczak, V. Volpert
Existence of reaction-diffusion-convection waves in unbounded
cylinders, Preprint, 2004.

\bibitem{B1} H.~Berestycki, The influence of advection on
the propagation of fronts in
reaction-diffusion equations, in
{\it Nonlinear PDEs in Condensed Matter and Reactive Flows},
NATO Science Series C, 569,
H. Berestycki and Y. Pomeau eds, Kluwer, Doordrecht, 2003.

\bibitem{BH} H.~Berestycki and F.~Hamel, Front propagation in
periodic excitable media,  Comm. Pure Appl. Math. {\bf 55}, 2002, 949--1032.

\bibitem{BHN-1} H. Berestycki, F. Hamel and N. Nadirashvili,
The speed of propagation for KPP type problems in periodic and general
domains, Preprint, 2003.

\bibitem{BHN-2} H. Berestycki, F. Hamel and N. Nadirashvili,
Elliptic eigenvalue problems with large drift and applications
to nonlinear propagation phenomena, Preprint, 2003.


\bibitem{BLN} H. Berestycki, B. Larrouturou and L. Nirenberg,
A nonlinear elliptic problem describing the propagation of a curved
premixed flame, in {\it Mathematical Modeling in Combustion and
Related Topics}, C.-M. Brauner and C. Schmidt-Lain\'e, eds., NATO ASI
Series, Kluwer, 1988.


\bibitem{BLL} H. Berestycki, B. Larrouturou and P. L. Lions,
 Multi-dimensional traveling wave solutions of a flame
propagation model, Arch. Rational Mech. Anal., {\bf 111},
1990, 33-49.

\bibitem{BNS}  H. Berestycki, B. Nicolaenco and B. Scheurer,
Traveling wave solutions to combustion models and their singular limits,
SIAM Jour. math. Anal., {\bf 16}, 1983, 1207-1242.

\bibitem{BN} H. Berestycki and L. Nirenberg,  Traveling fronts in
cylinders,  Annales de l'IHP, Analyse non lin\'eare, {\bf 9},
1992, 497-572.

\bibitem{CKOR} P. Constantin, A. Kiselev, A. Oberman, L. Ryzhik, Bulk
burning rate in passive-reactive diffusion, Arch. Rat. Mech. Anal. {\bf 154},
2000, 53-91.

\bibitem{CKR} P. Constantin, A. Kiselev and L. Ryzhik,
Fronts in reactive convection: bounds, stability and instability,
Comm. Pure Appl. Math., {\bf 56}, 2003, 1781-1803.

\bibitem{deWit} A. de Wit, Fingering of chemical fronts in porous media,
Phys. Rev. Let., {\bf 87}, 054502.


\bibitem{FKR} A. Fannjiang, A. Kiselev and L. Ryzhik, Unpublished notes, 2002.


\bibitem{Fisher} R. Fisher, The wave of advance of advantageous genes,
  Ann. Eugenics, {\bf 7}, 1937, 355--369.

\bibitem{F1} M. Freidlin and J.~G\"artner,  On the propagation
of concentration waves in periodic and random media, Soviet
Math. Dokl., {\bf 20}, 1979, 1282-1286.

\bibitem{F2} M.~Freidlin,  Geometric optics approach to
reaction-diffusion equations,  SIAM J. Appl. Math., {\bf 46},
1986, 222-232.

\bibitem{GT} D. Gilbarg and. Trudinger, Elliptic Partial Differential Equations
of Second Order, Springer-Verlag, 1983.

\bibitem{Hamel} F. Hamel, Formules min-max pour les vitesses d'ondes progressives multidimensionnelles,
Ann. Fac. Sci. Toulouse Math., S\'erie 6, {\bf 8}, 1999, 259--280.

\bibitem{HPS} S.~Heinze, G.~Papanicolau and A.~Stevens,
Variational principles for propagation speeds in inhomogeneous
media,  SIAM J. Appl. Math. {\bf 62}, 2001, 129-148.

\bibitem{KS} L.~Kagan and G.~Sivashinsky, Flame propagation and extinction
in large-scale vortical flows, Combust. Flame {\bf 120}, 2000, 222-232.

\bibitem{KRS}  L.~Kagan, P.D.~Ronney and G.~Sivashinsky,
Activation energy effect on flame propagation in large-scale vortical
flows, \rm Combust.  Theory Modelling {\bf 6}, 2002, 479-485.

\bibitem{KR} A.~Kiselev and L.~Ryzhik, {Enhancement of the travelling front
speeds in reaction-diffusion equations with advection}, Ann. Inst.
H. Poincar\'e Anal. Non Lin\'eaire {\bf 18}, 2001, 309-358.


\bibitem{KPP}
A.N.~Kolmogorov, I.G.~Petrovskii and N.S.~Piskunov,
\'Etude de l'\'equation de la chaleurde mati\`ere et son application \`a
un probl\`eme biologique, Bull. Moskov. Gos. Univ. Mat. Mekh. {\bf 1}
(1937), 1-25. (see \cite{curved-fronts} pp. 105-130 for an English transl.)



\bibitem{MS}
A.~Majda and P.~Souganidis, Large scale front dynamics for turbulent
reaction-diffusion equations with separated velocity scales,
Nonlinearity, {\bf 7}, 1994, 1-30.

\bibitem{MX} S. Malham and J. Xin, Global solutions to a reactive Boussinesq
system with front data on an infinite domain, Comm. Math. Phys., {\bf 193},
1998, 287-316.

\bibitem{Nash} J. Nash, Continuity of solutions of parabolic and
elliptic equations, Amer. Jour. Math., {\bf 80}, 1958, 931-954.

\bibitem{curved-fronts} Dynamics of curved fronts, P. Pelc\'e, Ed.,
  Academic Press, 1988.


\bibitem{PX-91}
G. Papanicolaou and X. Xin,
Reaction diffusion fronts in periodically layered media,  Jour. Stat. Phys.,
{\bf 63}, 1991, 915-932.



\bibitem{TV1} R.~Texier-Picard and V.~Volpert, Probl\`emes de
r\'eaction-diffusion-convection dans des cylindres non born\'es,
C. R. Acad. Sci. Paris Sér. I Math. {\bf 333}, 2001, 1077-1082

\bibitem{TV2} R.~Texier-Picard and V.~Volpert,
Reaction-diffusion-convection problems in unbounded cylinders,
Revista Matematica Complutense, {\bf 16}, 2003, ...

\bibitem{VR} N. Vladimirova, R. Rosner, Model flames in the Boussinesq limit:
the effects of feedback, Phys. Rev. E., {\bf 67}, 2003, 066305.

\bibitem{Volpert-1}V.A. Volpert and A.I. Volpert, {\it Location of spectrum
  and stability of solutions for monotone parabolic system}, Adv. Diff.
  Eq., {\bf 2}, 1997, 811-830.

\bibitem{Volpert-2}V.A. Volpert and A.I. Volpert, {\it Existence and
    stability of multidimensional travelling waves in the monostable
    case}, Israel Jour. Math., {\bf 110}, 1999, 269-292.

\bibitem{Volpert-3}V.A. Volpert and A.I. Volpert, {\it Spectrum of
    elliptic operators and stability of travelling waves}, Asymptotic
  Anal., {\bf 23}, 2000, 111-134.



\bibitem{Winn} B. Win, Ph.D. thesis, University of Chicago, 2004.

\bibitem{X1} J.~Xin, Existence of planar flame fronts in
convective-diffusive periodic media, Arch. Rat. Mech. Anal., {\bf
121}, 1992, 205-233.

\bibitem{X2} J.~Xin, Existence and nonexistence of travelling waves
  and reaction-diffusion front propagation in periodic media,
  Jour.  Stat. Phys., {\bf 73}, 1993, 893-926.

\bibitem{X3} J.~Xin, Analysis and modelling of front propagation in
heterogeneous media, SIAM Rev., {\bf 42}, 2000, 161-230.

\bibitem{ZBLM} Ya.B.~Zeldovich, G.I.~Barenblatt,
V.B.~Librovich and G.M.~Makhviladze,
\it The Mathematical Theory of Combustion and Explosions, \rm
Translated from the Russian by Donald H. McNeill. Consultants
Bureau [Plenum], New York, 1985.

\end{thebibliography}
\end{document}